%% file: 2001-7.tex
\documentclass{gtart}

\input gtoutput

\volumenumber{5}\papernumber{7}\volumeyear{2001}
\pagenumbers{227}{266}
\accepted{19 March 2001}
\proposed{Jean-Pierre Otal}
\seconded{Yasha Eliashberg, David Gabai}
\received{23 January 1999}
\revised{3 April 2000}
\published{23 March 2001}

\usepackage{graphicx,labelfig}

\newcommand{\mr}{\mbox{$\widetilde{M}\/$ }}

\newcommand{\fr}{\mbox{$\widetilde{\cal F}$ }}
\newcommand{\preua}{\noindent \mbox{\bf Proof\qua}}
\newcommand{\preuda}[1]
{\noindent \mbox{\bf Proof of \ref{#1}\qua}}

\newcommand{\fin}{\endproof}

\newtheorem{propa}{Proposition}[section]
\newtheorem{thma}[propa]{Theorem}
\newtheorem{lema}[propa]{Lemma}
\newtheorem{cora}[propa]{Corollary}

\theoremstyle{remark}
\newtheorem{defina}[propa]{Definition}
\newtheorem{raca}[propa]{Remark}

\newcommand{\rquea}{\begin{raca}}
\newcommand{\erquea}{\end{raca}}

\begin{document}  

\title{Flag Structures on Seifert Manifolds}

\author{Thierry Barbot}

\address{C.N.R.S, ENS Lyon, UMPA, UMR 5669
\\46, all\'ee d'Italie, 69364 Lyon, France}
\asciiaddress{C.N.R.S, ENS Lyon, UMPA, UMR 5669
\\46, allee d'Italie, 69364 Lyon, France}
\email{barbot@umpa.ens-lyon.fr}

\begin{abstract}
We consider faithful projective actions of a cocompact lattice of
$SL(2,{\bf R})$ on the projective plane, with the following property:
there is a common fixed point, which is a saddle fixed point for every
element of infinite order of the the group. Typical examples of such
an action are linear actions, ie, when the action arises from a
morphism of the group into $GL(2,{\bf R})$, viewed as the group of
linear transformations of a copy of the affine plane in ${\bf
R}P^{2}$.  We prove that in the general situation, such an action is
always topologically linearisable, and that the linearisation is
Lipschitz if and only if it is projective.  This result is obtained
through the study of a certain family of flag structures on Seifert
manifolds. As a corollary, we deduce some dynamical properties of the
transversely affine flows obtained by deformations of horocyclic
flows. In particular, these flows are not minimal.
\end{abstract}

\asciiabstract{We consider faithful projective actions of a
cocompact lattice of SL(2,R) on the projective plane, with the
following property: there is a common fixed point, which is a saddle
fixed point for every element of infinite order of the the
group. Typical examples of such an action are linear actions,
ie, when the action arises from a morphism of the group into GL(2,R),
viewed as the group of linear transformations of a copy of the affine
plane in RP^{2}.  We prove that in the general situation, such an
action is always topologically linearisable, and that the
linearisation is Lipschitz if and only if it is projective.  This
result is obtained through the study of a certain family of flag
structures on Seifert manifolds. As a corollary, we deduce some
dynamical properties of the transversely affine flows obtained by
deformations of horocyclic flows. In particular, these flows are not
minimal.}


\maketitlepage

\section{Introduction}

Let $\bar{\Gamma}$ be the fundamental group of a closed surface
with negative Euler characteristic. It admits many interesting
actions on the sphere $S^{2}$:

-\qua conformal actions through morphisms $\bar{\Gamma} \rightarrow PSL(2,{\bf C})$,

-\qua projective actions on the sphere of half-directions in ${\bf R}^{3}$
through morphisms $\bar{\Gamma} \rightarrow GL(3,{\bf R})$.

We have one natural family of morphisms from $\bar{\Gamma}$ 
into $PSL(2,{\bf C})$, and two natural families of morphisms from $\bar{\Gamma}$ 
into $GL(3,{\bf R})$:

\begin{enumerate}
\item {\em Fuchsian morphisms:}
fuchsian morphisms are faithful morphisms from
$\bar{\Gamma}$ into $PSL(2,{\bf R}) \subset PSL(2,{\bf C})$,
with image a cocompact discrete subgroup
of $PSL(2,{\bf R})$. In this case, the domain of discontinuity
of the corresponding action of $\bar{\Gamma}$
is the union of two discs, and these two discs
have the same boundary, which is nothing but the natural embedding
of the boundary of the Poincar\'e disc ${\bf H}^{2}$ into the boundary of the
hyperbolic $3$--space ${\bf H}^{3}$.
Moreover, on every component of the domain of discontinuity, 
the action of $\bar{\Gamma}$
is topologically conjugate to the action by isometries
through $PSL(2,{\bf R})$ on the Poincar\'e disc (the
topological conjugacy is actually quasi-conformal)
and the action on the common boundary of these discs is conjugate to
the natural action of $\bar{\Gamma}$ through $PSL(2,{\bf R})$ on
the projective line ${\bf R}P^{1}$. Finally, all these
actions on the whole sphere through $PSL(2,{\bf R})$ are 
quasi-conformally conjugate one to the other.

\item {\em Lorentzian morphisms:}
a lorentzian morphism is a faithful
morphism $\bar{\Gamma} \rightarrow SO_{0}(2,1) \subset GL(3,{\bf R})$
whose image is a cocompact lattice of $SO_{0}(2,1)$,
the group of linear transformations of determinant $1$ preserving
the Lor\-entzian cone of ${\bf R}^{3}$. Observe that such a
morphism corresponds to a fuchsian morphism {via} the isomorphism
$SO_{0}(2,1) \approx PSL(2,{\bf R})$.
The action on the projective plane associated to a lorentzian morphism has
the following properties:

-\qua it preserves an ellipse, on which the restricted action
is conjugate to the projective action
on ${\bf R}P^{1}$ through the associated fuchsian morphism,

-\qua it preserves a disc, whose boundary is the $\bar{\Gamma}$--invariant ellipse.
This disc is actually the projective Klein model of the Poincar\'e disc,
the action of $\bar{\Gamma}$ on it is conjugate to the associated fuchsian
action of $\bar{\Gamma}$ on the Poincar\'e disc,

-\qua it preserves a M\"{o}bius band (the complement of the closure of the invariant
disc). The action on it is topologically transitive (ie, there is
a dense $\bar{\Gamma}$--orbit). We have no need here to describe further 
this nice action.

Moreover, all the lorentzian actions are topologically conjugate 
one to the other,
and the conjugacy is H\"older continuous (we won't give any justification
here of this assertion, since it requires developments which are far away
from the real topic of this paper).

\item {\em Special linear morphisms:}
they are the faithful morphisms
$\bar{\Gamma} \rightarrow SL(2,{\bf R})$, where
$SL(2,{\bf R})$ is considered here as
the group $SL \subset SL(3,{\bf R})$ of matrices of positive determinant and of the form:
\[ \left(\begin{array}{ccc}
\ast & \ast & 0 \\
\ast & \ast & 0 \\
0  & 0 & 1
\end{array}\right)\]
Moreover we require that the image of the morphism is a lattice in $SL$.
Then, the action of $\bar{\Gamma}$ on the projective plane 
has a common fixed point, an invariant projective line, and an 
invariant punctured affine plane.
The action on the invariant line is the usual projective action
on ${\bf R}P^{1}$ through the natural projection $SL \rightarrow PSL(2,{\bf R})$,
and the action on the punctured affine plane
is the usual linear action. This action is minimal (every
orbit is dense) and uniquely ergodic (there is an unique invariant
measure up to constant factors). Contrary to the preceding
cases, the action highly depends on the morphism into 
$SL$: two morphisms induce topologically conjugate
actions if and only if they are conjugate by an inner automorphism
in the target $SL$. 

\end{enumerate}

We are interested in the small deformations of these actions
arising by perturbations of the morphisms into $PSL(2,{\bf C})$
or $GL(3,{\bf R})$. We list below the main properties of
these deformed actions; we will see later how to justify all these
claims.

\begin{enumerate}

\item {\em Quasi-fuchsian actions:\/} morphisms from
$\bar{\Gamma}$ into $PSL(2,{\bf C})$ which are
small deformations of fuchsian morphisms are {\em quasi-fuchsian\/}: 
this essentially means
that their associated actions on the sphere 
are quasi-conformally conjugate to fuchsian actions.
They all preserve a Jordan curve, this Jordan curve is rectifiable if
and only if it is a great circle, in which case the action
is actually fuchsian (see for example \cite{sulli}, chapter $7$).

\item {\em Convex projective actions:\/} 
we mean by this the actions arising from morphisms from
$\bar{\Gamma}$ into $GL(3,{\bf R})$ near 
lorentzian projective morphisms. Such an action
still preserves a strictly convex subset of ${\bf R}P^{2}$ whose
boundary is a Jordan curve of class $C^{1}$ (it is of class
$C^{2}$ if and only if the action is conjugate in $PGL(3,{\bf R})$
to a lorentzian action, see \cite{benzecri}). Moreover, all these actions 
are still
topologically conjugate one to the other\footnote{In this case, we have the additional 
remarkable fact: in the variety of morphisms $\bar{\Gamma} \rightarrow 
PGL(3,{\bf R})$, the morphisms belonging to the whole connected component of
the lorentzian morphisms (the so-called Hitchin component) 
induce the same action on the projective
plane up to topogical conjugacy \cite{choigold}.}.

\item {\em Hyperbolic actions:\/} these are the real topic of this paper,
thus we discuss them below in more detail.

\end{enumerate}

Hyperbolic actions arise from morphisms
from $\bar{\Gamma}$ into $PGL(3,{\bf R})$ which are 
deformations of special linear morphisms.
Actually, we will not consider all these deformations;
we will restrict ourselves to the deformations for which the deformed
action has still an invariant point: they correspond to morphisms
into the group $Af^{\ast}_{0}$ of matrices of the form:
\[ \left(\begin{array}{cc}
A &\begin{array}{c} 
0 \\
0
\end{array}\\
\begin{array}{cc}
x & y 
\end{array} & 1
\end{array}\right)\]
where $A$ is a $2 \times 2$--matrix of positive determinant 
(we will say that the matrix $A$ is the 
linear part, and that $(x,y)$ is the translation part).
This group is in a natural way dual to the group $Af_{0}$ of 
orientation preserving
affine transformations of the plane: the space of projective 
lines in ${\bf R}P^{2}$
is a projective plane too, and the dual action of $Af^{\ast}_{0}$
on this dual projective plane preserves a projective copy of the affine
plane.

Small deformations $\bar{\Gamma} \rightarrow Af_{0}^{\ast}$ of special linear morphisms 
all satisfy the following properties (cf
Lemma \ref{ouvR}):

-\qua the morphism $\bar{\Gamma} \rightarrow Af^{\ast}_{0}$ is injective,

-\qua the common fixed point is a fixed point of saddle type
for every non-trivial element of $\bar{\Gamma}$.
Equivalently, the image of every non-trivial element of
$\bar{\Gamma}$ in the dual group $Af_{0}$ is a hyperbolic
affine transformation.

Morphisms $\rho\co  \bar{\Gamma} \rightarrow Af^{\ast}_{0}$ satisfying the
properties above are called {\em hyperbolic.\/} 
In the special case where the translation part $(x,y)$ is zero
for every element, we say that the hyperbolic action is {\em horocyclic\/}
(we will soon justify this terminology).
Observe that the conjugacy by homotheties of the form
\[ \left(\begin{array}{ccc}
e^{t} & 0 & 0 \\
0 & e^{t} & 0 \\
0 & 0 & 1
\end{array}\right) \]
does not modify the linear parts, but multiply the translation part
$(x,y)$ by $e^{t}$. It follows that hyperbolic morphisms can all
be considered as small deformations of horocyclic morphisms (cf Proposition
\ref{fuch}).

Hyperbolic morphisms can
be defined in another way: we call the {\em unimodular linear part of $\rho$\/}
the projection in $SL(2,{\bf R})$ of the linear part of the morphism;
we denote it by $\rho_{0}$. For every element $\gamma$ of $\bar{\Gamma}$,
let $\bar{u}(\gamma)$ be the logarithm of the determinant of the linear part 
of $\rho(\gamma)$
(as an linear transformation of the plane). 
It induces an element of $H^{1}(\bar{\Gamma},{\bf R})$. On the other hand,
$H^{1}(\bar{\Gamma}, {\bf R})$ is isomorphic to $H^{1}(\Sigma, {\bf R})$,
where $\Sigma$ is the quotient of the Poincar\'e disc by the projection
of $\rho_{0}(\bar{\Gamma})$ in $PSL(2,{\bf R})$. The surface $\Sigma$
is naturally equipped with a hyperbolic metric, and thus, 
we can consider the {\em stable norm on $H^{1}(\Sigma, {\bf R})$\/}
(this stable norm depends on $\rho_{0}$)
Then (Remark \ref{defhyp}),
the morphism $\rho$ is hyperbolic if and only if the morphism $\rho_{0}$
is fuchsian (ie, has a dicrete cocompact image),
and if the stable norm of $\bar{u}$ is less than $\frac{1}{2}$.
We call hyperbolic every projective action of $\bar{\Gamma}$ induced
by a hyperbolic morphism.
The main result of this paper is (Corollaries \ref{conjug}, \ref{lip}):

\medskip\noindent\mbox{\bf Theorem A}\qua
\sl
Every hyperbolic action of $\bar{\Gamma}$ is topologically
conjugate to the projective horocyclic action of its linear part.
The conjugacy is Lipschitz if and only if it is a projective transformation.
\rm\medskip

As a corollary, any hyperbolic action preserves an annulus on
which it is uniquely ergodic, and the two boundary components
of this annulus
are respectively the common fixed point and an invariant Jordan curve
(Corollary \ref{Uh}). 
We give below a computed picture of such a Jordan curve:

\begin{figure}[htb]
\label{lambda}
 \centerline{\includegraphics[width=5cm, height=5cm]{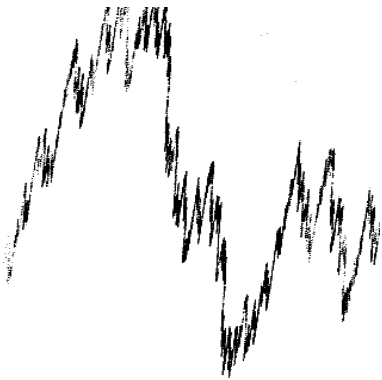}}
\caption{A zoom on the invariant Jordan curve}
\end{figure}

The studies of all these deformations have a common feature:
we have to transpose the problem to a $3$--dimensional object.

\begin{enumerate}

\item {\em The case of fuchsian actions:\/} in this case, the 
key idea is to consider
the quotient of hyperbolic $3$--space ${\bf H}^{3}$ by
$\bar{\Gamma}$ (viewed as a subgroup of $PSL(2,{\bf C}) \approx Isom({\bf H}^{3})$).
It is a hyperbolic $3$--manifold, homeomorphic to the product of
a surface $\Sigma$ by $]0,1[$. The action in ${\bf H}^{3}$ has
a finite fundamental polyhedron (see \cite{marden}, chapter $4$).
The fuchsian morphism can be considered as the holonomy morphism
of this hyperbolic manifold. 
It is well-known that any deformation of the holonomy corresponds to
a deformation of the hyperbolic structure (this is
a general fact about $(G,X)$--structures,
see for example \cite{golbou}, \cite{canary}). 
According to \cite{marden}, Theorem $10.1$, the deformed action
still has a finite sided polyhedron. It follows then that the domain
of discontinuity of the deformed action contains two invariant discs,
and then, that the action is quasi-fuchsian, ie, that it is quasi-conformally
conjugate to a fuchsian action (\cite{marden}, section $3.2$).
(Quasi-conformal stability of quasi-fuchsian groups
is also proved by L Bers in \cite{bers}, using
different tools). 

\item {\em The case of convex projective actions:\/}
the deformations of lorentzian cones can be understood by the
following method: the invariant disc is the
projection in ${\bf R}P^{2}$ of the lorentzian cone.
Add to the cocompact lattice in $SO_{0}(2,1)$ any homothety
of ${\bf R}^{3}$ of non-constant factor. 
We obtain a new group which acts freely, properly and cocompactly
on the lorentzian cone. The quotient of this action 
is a closed $3$--manifold, equipped with a 
{\em radiant affine structure,\/} ie, a $(GL(3,{\bf R}),{\bf R}^{3})$--structure.
It follows from a Theorem of J\,L Koszul \cite{kosz} that for
any deformation of the holonomy morphism, the corresponding
deformed radiant affine manifold is still the quotient of some
convex open cone in ${\bf R}^{3}$. It provides the invariant
strictly convex subset in ${\bf R}P^{2}$. We won't discuss here
why the $\bar{\Gamma}$--action is still conjugate to the
lorentzian action.

\item {\em The case of hyperbolic actions:\/}
we will deal with this case by considering 
{\em flag manifolds.\/}

\end{enumerate}

A flag manifold is a closed $3$--manifold equipped with a $(G,X)$--structure
where the model space $X$ is the flag variety, ie, the
set of pairs $(x,d)$, where $x$ is a point of the projective
plane, and $d$ is an oriented projective line through $x$. The group 
$G$ to be considered is the group $PGL(3,{\bf R})$
of projective transformations.
A typical example of such a structure is given
by the projectivisation of the tangent bundle
of a $2$--dimensional real projective orbifold. This family
is fairly well-understood, thanks to the
classification of compact real projective surfaces
(see \cite{cdcr1, cdcr2, cdcr3, CG}). 
Anyway,
the flag manifolds we will consider here are
of different nature.

The prototypes of the flag manifolds we will consider here are obtained in
the following way: 
consider the $GL_{0}$--invariant copy of the affine plane ${\bf R}^{2}$
in ${\bf R}P^{2}$, and let $0$ be the fixed point of $GL_{0}$ in ${\bf R}^{2}$. Let
$X_{0}$ be the open subset of $X$ formed by the pairs
$(x,d)$, where $x$ belongs to ${\bf R}^{2} \setminus \{ 0 \}$,
and $d$ is a projective line containing $x$ but not $0$.
Then, the subgroup $SL(2,{\bf R}) \subset Af_{0}^{\ast} \subset PGL(3,{\bf R})$
acts simply transitively on $X_{0}$. Therefore, 
if $\rho_{0}\co  \bar{\Gamma} \rightarrow SL(2,{\bf R})$ is a faithful
morphism with discrete and cocompact image,
the $\bar{\Gamma}$--action on $X_{0}$ through $\rho_{0}$ is free
and properly discontinuous. The quotient of this action
is a flag manifold, homeomorphic to the unitary tangent bundle of a surface.
Actually, it follows from a Theorem of F Salein that 
horocyclic actions on $X_{0}$ are free and properly discontinuous too
(Corollary \ref{sal}).
We call {\em canonical Goldman flag manifolds\/} all the quotient
manifolds of actions obtained in this way.
In this case, the morphism $\rho_{0}$ is not strictly speaking the
holonomy morphism of the flag structure, because $\bar{\Gamma}$
is not the fundamental group $\Gamma$ of the flag manifold,
but the quotient of it by its center. We will actually
consider the morphism $\Gamma \rightarrow GL_{0}$ induced
by $\rho_{0}$; and we will still denote it by $\rho_{0}$.
Then, the definition of hyperbolic morphism has to
be generalised for morphisms $\Gamma \rightarrow Af^{\ast}_{0}$
(cf section \ref{notation}).

By deforming the morphism $\rho_{0}$, we obtain new flag manifolds.
Small deformations still satisfy:

-\qua the ambient flag manifold is homeomorphic to the unitary tangent
bundle of a surface,

-\qua the holonomy morphism is hyperbolic,

-\qua the image of the developing map is contained in $X_{\infty}$,
the open subspace of $X$ formed by the pairs $(x,d)$ where
$x$ belongs to ${\bf R}P^{2} \setminus \{ 0 \}$ and where
$d$ does not contain $0$ (see section \ref{flag}).

We call flag manifolds satisfying these $3$ properties {\em Goldman
flag manifolds\/} .  The main step for the proof of Theorem A is the
following theorem (section \ref{flagol}):

\medskip\noindent\mbox{\bf Theorem B}\qua
\sl
Let $M$ be a Goldman flag manifold with holonomy morphism $\rho$.
Then, $M$ is the quotient of an open subset $X(\rho)$
of $X_{\infty} \subset X$ which has the following description:
there is a Jordan curve $\Lambda(\rho)$ in ${\bf R}P^{2}$
which does not contain the common fixed point $0$, and
$X(\rho)$ is the set of pairs $(x,d)$ where $x$
belongs to ${\bf R}P^{2} \setminus (\Lambda(\rho) \cup \{ 0 \})$
and $d$ does not contain $0$.
\rm\medskip

Any flag manifold 
in\-he\-rits two $1$--dimensional foliations, that we call
the {\em tautological foliations.\/} They arise from
the $PGL(3,{\bf R})$--invariant tautological foliations on $X$ whose
leaves are the $(x,d)$ where $x$ and $d$ respectively remain fixed.
The tautological foliations are naturally transversely
real projective. We observe only in this introduction
that projectivisations of tangent bundles of real
projective orbifolds can be characterized
as the flag manifolds such that one of their
tautological foliations has only compact leaves
(this observation has
no incidence in the present work).  

In the case of canonical Goldman flag manifolds,
the tautological foliations are transversely affine.
Actually, they are the horocyclic foliations associated
to the exotic Anosov flows defined in \cite{ghyqf2}.
This justifies our
terminology ``horocyclic actions", the fact that
horocyclic actions are uniquely ergodic (since horocyclic
foliations of exotic Anosov flows are uniquely ergodic \cite{bowen}),
and the non-conjugacy between different horocyclic actions
(since horocyclic foliations are rigid (cf \cite{abe})). 

When the Goldman flag manifold is {\em pure,\/} ie, when
it is not isomorphic to a canonical flag Goldman manifold,
one of these foliations is no longer transversely affine; in fact
we understand this foliation quite well, since 
it is topologically conjugate to an exotic horocyclic foliation
(Theorem \ref{psi}).

The situation is different for the other tautological foliation:
they have been first introduced by W Goldman, which defined
them as the flows obtained by deformation of
horocyclic foliations amongst transversely affine foliations
on a given Seifert manifold $M$ (the two definitions coincide, 
see Proposition \ref{tra} and the following discussion).
For this reason, we call these foliations {\em Goldman foliations,\/}
and we extend this terminology to the ambient flag manifold.
As observed by S Matsumoto \cite{matsu}, 
nothing is known about the dynamical properties of pure Goldman foliations, even
when they preserve a transverse parallel volume form.
As a consequence of this work, we can prove
(section \ref{sixdeux}):

\medskip\noindent\mbox{\bf Theorem C}\qua
\sl
Goldman foliations are not minimal.
\rm\medskip

Hence, the dynamical properties of
Goldman foliations are drastically different from
the dynamical properties of horocyclic foliations.

Finally, many questions on the subject are still open.
The presentation of these problems is the conent of the last section
(Conclusion) of this paper.

Special thanks are due to to Damien Gaboriau, Jean-Pierre Otal and
Abdelghani Zeghib for their valuable help.
 
\section{Preliminaries}

\subsection{Notation}
\label{notation}
$M$ is an oriented 
closed $3$--manifold. We denote by $p\co  \widetilde{M} \rightarrow M$
a universal covering and $\Gamma$ the Galois group of this covering, ie,
the fundamental group of $M$. 

We denote by ${\bf R}P^{2}$ the usual projective plane, and
${\bf R}P^{2}_{\ast}$ its dual: ${\bf R}P^{2}_{\ast}$ is the
set of projective lines in ${\bf R}P^{2}$.
Let $\kappa\co  {\bf R}P^{2} \rightarrow {\bf R}P^{2}_{\ast}$ be the
duality map induced by the identification of ${\bf R}^{3}$
with its own dual, mapping the canonical basis of ${\bf R}^{3}$
to its canonical dual base.
Since ${\bf R}^{3}$ is also the dual space of its own dual, we obtain by
the same way an isomorphism $\kappa^{\ast}\co  {\bf R}P^{2}_{\ast} \rightarrow {\bf R}P^{2}$, which is the inverse of $\kappa$.

We denote by $X$ the flag variety: this is the subset of 
${\bf R}P^{2} \times {\bf R}P^{2}_{\ast}$
formed by the pairs $(x,d)$ where $d$ is an projective line containing $x$.
Let $p_{1}$ and $p_{2}$ be the projections of $X$ over ${\bf R}P^{2}$ and 
${\bf R}P^{2}_{\ast}$.
The flag variety $X$ is naturally identified with the 
projectivisation of the tangent
bundle of ${\bf R}P^{2}$. 
Let $\Theta$ be the orientation preserving involution of $X$ 
defined by $\Theta(x,d)=(\kappa^{\ast}(d), \kappa(x))$.

Let $PGL(3,{\bf R})$ be the group of 
projective automorphisms
of ${\bf R}P^{2}$. 
The differential of the action of $PGL(3,{\bf R})$ on ${\bf R}P^{2}$ induces an
orientation preserving action on $X$.
Consider the Cartan involution on $GL(3,{\bf R})$ mapping a matrix 
to the inverse of its transposed matrix. It induces an involution
$\theta$ of $PGL(3,{\bf R})$.
We have the equivariance relation $\Theta \circ A = \theta(A) \circ \Theta$
for any element $A$ of $PGL(3,{\bf R})$.

A flag structure on $M$ is a $(PGL(3,{\bf R}), X)$--structure on $M$ in
the sense of \cite{ratcliffe}.  We denote by ${\cal D}\co  \widetilde{M}
\rightarrow X$ its developing map, and by $\rho\co  \Gamma \rightarrow
PGL(3,{\bf R})$ its holonomy morphism.  The compositions of $\cal D$
and $\rho$ by $\Theta$ and $\theta$ define another flag structure on
$M$: the {\em dual flag structure.\/} In general, a flag structure is
not isomorphic to its dual.

On $X$, we have two natural one dimensional foliations by circles:
the foliations whose leaves are the fibers of $p_{1}$ and $p_{2}$.
We call them respectively the {\em first\/}
and
the {\em second tautological foliation.\/} 
They are both preserved by the action of $PGL(3,{\bf R})$.
Therefore, they induce on each manifold equipped
with a flag structure two foliations that we still
call the first and second
tautological foliations.
The first (respectively second) tautological foliation
is the second (respectively first) tautological foliation
of the dual flag structure. Observe that these
foliations are transversely real projective.
Observe also that they are nowhere
collinear, and that the plane field that
contains both is a contact plane field.

Consider the usual embedding of the affine plane ${\bf R}^{2}$ in $P^{2}{\bf R}$.
We denote by $0$ the origin
of ${\bf R}^{2}$.
The boundary of ${\bf R}^{2}$ in ${\bf R}P^{2}$ is the projective line $\kappa(0)$, the
line at infinity. We denote it by $d_{\infty}$. It is naturally identified
with the set ${\bf R}P^{1}$ of lines in ${\bf R}^{2}$ through $0$.
We identify thus the group of transformations
of the plane with the group of projective
transformations preserving the line $d_{\infty}$. 
Let $Af_{0}$ be the group of orientation preserving affine transformations.
The elements of $Af_{0}$ are the projections in $PGL(3,{\bf R})$
of matrices of the form:
\[
\left(\begin{array}{cc}
A & \begin{array}{c}
     u \\
     v
    \end{array}\\
\begin{array}{cc}
0 & 0 
\end{array} & 1
\end{array}\right)
\]
where $A$ belongs $GL_{0}$, the group of $2 \times 2$ matrix with positive determinant.
The group 
$GL_{0}$ is the stabilizer in $Af_{0}$ of the point
$0$. We denote by $SL$ the subgroup formed by the
elements of $GL_{0}$ of determinant $1$ (as a group of linear transformation
of the plane; equivalently, $SL$ is the derived subgroup of $GL_{0}$), 
by $p_{0}\co  \widetilde{SL} \rightarrow SL$ 
the universal covering map, and by $PSL$ the quotient of
$SL$ by its center $\{ \pm Id \}$.

Let $\Gamma$ be a cocompact lattice of $\widetilde{SL}$.
Let $H$ be the center of $\Gamma$. We select a generator $h$ of $H \approx \bf Z$.
Let $\bar{\Gamma}$ be the quotient of $\Gamma$ by $H$.
We denote by $\rho_{0}\co  \Gamma \rightarrow \bar{\Gamma} \subset SL \subset Af_{0}$
the quotient map: this is the restriction of $p_{0}$ to $\Gamma$.

Let ${\cal R}({\Gamma})$ be the space 
of representations of ${\Gamma}$
into $Af_{0}$.
It has a natural structure of an algebraic variety.

Let $Rep(\Gamma, PSL)$ be the space of morphisms of
$\Gamma$ into $PSL$. The elements of $Rep(\Gamma, PSL)$
vanishing on $H$ form a subspace that we denote by $Rep(\bar{\Gamma}, PSL)$.
As suggested by the notation, $Rep(\bar{\Gamma}, PSL)$ can be identified with the space
of representations of $\bar{\Gamma}$ into $PSL$. 

By taking the linear part of $\rho(\gamma)$, and then projecting in
$PGL_{0} \approx PSL$, we define an open map 
$\lambda\co  {\cal R}({\Gamma}) \rightarrow Rep({\Gamma},PSL)$.
We call $\lambda(\rho)$ the {\em projectivised linear part of $\rho$.\/}

An element $\rho$ of ${\cal R}({\Gamma})$
is {\em hyperbolic\/} if it satisfies the following conditions:

-\qua the kernel of $\lambda(\rho)$ is $H$,

-\qua for every element $\gamma$ of $\Gamma$ which has no non-trivial power belonging to $H$,
$\rho(\gamma)$ has two real eigenvalues,
one of absolute value strictly greater than $1$, and the other of absolute value 
strictly less than $1$. In other words, $\rho(\gamma)$ has a fixed point 
of saddle type.

Observe that this definition is dual to the definition given in the
introduction. 
A typical example of hyperbolic representations is $\rho_{0}$.
We denote by ${\cal R}_{h}({\Gamma})$ 
the set of 
elements of ${\cal R}({\Gamma})$ which are hyperbolic.

Let ${\cal T}(\bar{\Gamma})$ be the space of {\em cocompact fuchsian representations of
$\bar{\Gamma}$ into $PSL$\/}, ie, injective representations with
a discrete and cocompact image in $PSL$.
It is well-known that it is a connected component of the
space $Rep(\bar{\Gamma}, PSL)$ of all representations 
$\bar{\Gamma} \rightarrow PSL$.

\begin{lema}
\label{ouvR}
${\cal R}_{h}(\Gamma)$ is
an open subset of
${\cal R}({\Gamma})$. Its image by $\lambda$ is ${\cal T}(\bar{\Gamma})$.
\end{lema}

\preua
Let $Rep_{0}(\Gamma, PSL)$ be the subspace of $Rep(\Gamma, PSL)$ formed by the
morphisms $\rho$ with non-abelian image. 
This is an open subspace. For any element $\rho$ of  $Rep(\Gamma, PSL)$,
the image of $\rho$ is contained in the centralizer
of $\rho(h)$. But
the centralizers of non-trivial elements of $PSL$ are
all abelian, thus $Rep_{0}(\Gamma, PSL)$ is an open
subset of $Rep(\bar{\Gamma}, PSL)$.
Moreover, $Rep_{0}(\Gamma, PSL)$ obviously contains ${\cal T}(\bar{\Gamma})$.

Take any element $\rho$ of ${\cal R}_{h}({\Gamma})$.
Since $\bar{\Gamma}$ is not abelian, and since 
the kernel of $\rho$ is contained in $H$,
$\lambda(\rho)$
belongs to $Rep(\bar{\Gamma}, PSL)$. Moreover, 
$\lambda(\rho)\co  \bar{\Gamma} \rightarrow PSL$ is injective.
Let $N_{0}$ be the identity component
of the closure of $\lambda(\rho)(\bar{\Gamma})$ in $PSL$.
Then, $\lambda(\rho)^{-1}(N_{0} \cap \lambda(\rho)(\bar{\Gamma}))$ is 
a normal subgroup of $\bar{\Gamma}$. Hence, either it is
contained in the center $H$, or it is not solvable.
In the second case, $N_{0}$ is not solvable too: it must contain
elliptic elements with arbitrarly small rotation angle. 
But $\rho({\Gamma})$ contains then many elliptic elements
with rotation angles arbitrarly small: this is a contradiction since $\rho$
is hyperbolic.

Therefore, $\lambda(\rho)^{-1}(N_{0} \cap \rho(\bar{\Gamma}))$
is trivial, ie, $\rho(\Gamma)$ is discrete. Since $\lambda(\rho)(\bar{\Gamma})$
is isomorphic to $\bar{\Gamma}$, its cohomological dimension is
two. Hence, it is a cocompact subgroup of $PSL$, and
$\lambda({\cal R}_{h}(\bar{\Gamma}))$
is contained in ${\cal T}(\bar{\Gamma})$. 
The lemma follows. \fin

\rquea
\label{defhyp}
Lemma \ref{ouvR} enables us to give a method for defining all
hyperbolic morphisms: take any cocompact fuchsian
group $\bar{\Gamma}$ in $PSL$, and let $\tilde{\Gamma}$ be
the preimage by $p_{0}$ of $\bar{\Gamma}$. Take any
finite index subgroup $\Gamma$ of $\tilde{\Gamma}$. 
Denote by $\rho_{0}$ the restriction of $p_{0}$ to $\Gamma$.
Take now any morphism $u$ from $\Gamma$ into the multiplicative group
${\bf R} \setminus 0$.
We can now define a new morphism $\rho_{u}\co  \Gamma \rightarrow GL_{0}$
just by requiring
$\rho_{u}(\gamma) = u(\gamma) \rho_{0}(\gamma)$. Actually, all the
$\rho_{u}$ are nothing but the elements of the fiber of $\lambda$
containing $\rho_{0}$. 
The absolute value of the morphism $u$
is a morphism $\vert u \vert\co  {\Gamma} \rightarrow {\bf R}^{+}$.
Since $h$ admits a non-trivial power belonging to the commutator
subgroup $[\Gamma, \Gamma]$, $\vert u \vert$ is trivial on $H$;
therefore, it induces a morphism $\bar{u}\co  \bar{\Gamma} \rightarrow {\bf R}^{+}$.

Now, the following claim is easy to check:
{\sl the morphism $\rho_{u}$ is hyperbolic if and only if for any non-elliptic element
$\gamma$ of $\bar{\Gamma}$, the absolute value $\bar{u}(\gamma)$ belongs
to $]r(\gamma)^{-1}, r(\gamma)[$, where $r(\gamma)$ is the spectral radius of $\gamma$.\/}

This condition can be expressed in a more elegant way:
the logarithm of $\bar{u}$ is a morphism 
$L_{u}\co  \bar{\Gamma} \rightarrow {\bf R}$, ie, an element of 
$H^{1}(\bar{\Sigma}, {\bf R})$. On this cohomology space, we have
the {\em stable norm\/} (cf \cite{norm}) which is defined as
follows: for any hyperbolic element $\gamma$ of $\bar{\Gamma}$,
let $t(\gamma)$ be the double of the logarithm of $r(\gamma)$
(this is the length of the closed geodesic associated to $\bar{\Gamma}$
in the quotient of the Poincar\'e disc by $\bar{\Gamma}$).
For any element $\hat{\gamma}$ of $H_{1}(\bar{\Sigma}, {\bf Z})$,
and for any positive integer $n$, let $t_{n}(\hat{\gamma})$ the
infimum of the $\frac{t(\gamma)}{n}$ where $\gamma$ describes
all the elements of $\Gamma$ representing $n\hat{\gamma}$.
The limit of $t_{n}(\hat{\gamma})$ exists, it is the {\em stable
norm of $\hat{\gamma}$ in $H_{1}(\bar{\Sigma}, {\bf Z})$.\/} 
This norm is extended in an unique way on all 
$H_{1}(\bar{\Gamma}, {\bf R})$; the dual of it is the
{\em stable norm of $H^{1}(\bar{\Sigma}, {\bf R})$.\/}
The proof of the following claim is left to the reader:
{\sl the representation $\rho_{u}$ is hyperbolic if and only if 
the stable norm of $L_{u}$ is strictly less than $\frac{1}{2}$.\/}

\erquea

\rquea
\label{notorsion}
According to Selberg's Theorem, asserting that any finitely generated linear group
admits a finite index subgroup without torsion, for any hyperbolic
representation $\rho\co  \Gamma \rightarrow Af_{0}$, there exists a finite
index subgroup $\Gamma'$ of $\Gamma$ on which $\rho$ restricts as a hyperbolic
representation. This hyperbolic representation has the following properties:

-\qua its kernel is precisely the center of $\Gamma'$,

-\qua every non-trivial element of $\rho(\Gamma')$ is hyperbolic.

\erquea

Let $Af_{0}^{\ast}$ be the dual ${\theta}(Af_{0})$ of $Af_{0}$. 
Since $Af_{0}$ preserves the line at infinity $d_{\infty}$, the group
$Af_{0}^{\ast}$ fixes
the point $0$ in ${\bf R}P^{2}$.
It preserves also the open set $X_{\infty}$ whose
elements are the pairs $(x,d)$, where $x$ is a point
of ${\bf R}P^{2} \setminus 0$, and $d$ a line containing
$x$ but not $0$. Observe that the fundamental group of $X_{\infty}$
is infinite cyclic.
The group $GL_{0}$ (which is equal to its dual $\theta(GL_{0})$) preserves 
the subset
$X_{0} \subset X_{\infty}$
where $(x,d)$ belongs to $X_{0}$ if and only if $x$ belongs
to ${\bf R}^{2} \setminus \{ 0 \}$, and $d$ does not contain
$0$. Actually, the action of $SL$ on $X_{0}$ is simply transitive.
A representation $\rho\co  {\Gamma} \rightarrow Af_{0}^{\ast}$
is said to be hyperbolic if it is the dual representation of an
element of ${\cal R}_{h}(\Gamma)$.
Equivalently, it means that the point $0$ 
is a fixed point of saddle type of every $\rho(\gamma)$, when
$\gamma$ is of of infinite order.
Such a representation is given by a morphism $\rho_{1}\co  \Gamma \rightarrow GL_{0}$ and
two cocycles $u$ and $v$ such that $\rho(\gamma)$ is 
the projection in $PGL(3,{\bf R})$ of:
\[
\left(\begin{array}{cc}
\rho_{1}(\gamma) & \begin{array}{c}
     0 \\
     0
    \end{array}\\
\begin{array}{cc}
u(\gamma) & v(\gamma) 
\end{array} & 1
\end{array}\right)
\]
The morphism $\rho_{1}$ is the linear part of
$\rho$. It is a {\em horocyclic morphism.\/}

An $Af_{0}$--foliation is a foliation admitting a transverse 
$(Af_{0}, {\bf R}^{2})$--structure.

\subsection{Convex and non-convex sets}

Here, we collect some elementary facts on affine manifolds.

\begin{defina}
Let $X$ be a flat affine manifold. An open subset $U$ of $X$
is convex if any pair $(x,y)$ of points of $U$ are extremities
of some linear path contained in $U$. 
The exponential ${\cal E}_{x}$ of a point $x$ of $X$ is the
open subset of $X$ formed by the points which are extremities
of linear paths starting from $x$.
\end{defina}

The following lemmas are well-known. A good reference is \cite{carriere}.

\begin{lema}
\label{conv}
The developing map of a flat convex simply connected affine manifold
is a homeomorphism onto its image. \fin
\end{lema}

\begin{lema}
\label{expcon}
Let $X$ be a connected flat affine manifold.
If the exponential of every point of $X$ is convex,
then $X$ is convex. \fin
\end{lema}

\begin{lema}
\label{interconv}
Let $X$ be a flat affine simply connected manifold.
Let $U$ and $V$ be two convex subsets of $X$. If
$U \cap V$ is not empty, the restriction
of the developing map $\cal D$ to $U \cup V$ 
is a homeomorphism over ${\cal D}(U) \cup {\cal D}(V)$.\fin
\end{lema}

\begin{lema}
\label{conplan}
Let $U$ be an open 
star-shaped neighborhood of a point $x$ in the plane. If $U$ is not
convex, then it contains two points $y$ and $z$ such that:
\begin{itemize}
\item $x$, $y$ and $z$ are not collinear,
\item the closed triangle with vertices $x$, $y$ and $z$ is not contained
in $U$, 
\item the open triangle with vertices $x$, $y$ and $z$, and
the sides $[x,y]$, $[x,z]$, are contained in $U$.
\end{itemize}
\end{lema}
   
\preua
Let $y'$ and $z$ be two points of $U$ such that
the segment $[y',z]$ is not contained in $U$.
Observe that $x$, $y'$ and $z$ are not collinear.
Let $T_{0}$ be the closed triangle of vertices $x$, $y'$ and $z$.
For any real $t$ in the interval $[0,1]$, let $y_{t}$ be the point
$ty'+(1-t)x$. Let $I$ be the set of parameters
$t$ for which the segment $[y_{t},z]$ is contained in $U$.
It is open, non-empty since $0$ belongs to it, and 
does not contain $1$. Let $t$ be a boundary point of $I$:
the points $y_{t}$ and $z$ have the properties
required by the lemma.\fin

\section{Existence of flag structures}
\label{flag}

Let $M$ be a {\em principal Seifert manifold,\/} ie, 
the left quotient
of $\widetilde{SL}$ by a cocompact lattice $\Gamma$.
Let $\bar{\Gamma}$ be the projection of $p_{0}(\Gamma)$ in $PSL$.
Topologically, $M$ is a Seifert bundle over
the hyperbolic orbifold $\bar{\Sigma}$, quotient of the Poincar\'e disc
by $\bar{\Gamma}$.

Choose any element $v$ of $X_{0}$. Consider the map $\widetilde{SL}
\rightarrow X_{0} \subset X$ that maps $g$ to $p_{0}(g)(v)$, and the
morphism $\rho_{0}\co  \Gamma \rightarrow PGL(3,{\bf R})$, which is the
composition of $p_{0}$ with the inclusion $SL \subset Af_{0}^{\ast}
\subset PGL(3,{\bf R})$.  They are the developing map and holonomy
morphism of some flag structure on $M$. Observe that this structure
does not depend on the choice of $v$. We call the flag structures
obtained in this way the {\em unimodular canonical flag structures.\/}

We are concerned here with the deformations of unimodular
canonical flag structures.
Let $t \mapsto \rho_{t}$ be a deformation of $\rho_{0}$
inside $PGL(3,{\bf R})$, where the
parameter
$t$ belongs to $[0,1]$. 
As we recalled in the introduction, for small $t$, the
morphisms $\rho_{t}$ is the holonomy morphism of some new flag structure.
Moreover, these deformed flag structures
near the canonical one are well-defined up to isotopy
by their holonomy morphisms.
We are interested by the deformations of $\rho_{0}$ inside $Af_{0}^{\ast}$,
ie, where all the $\rho_{t}$ are morphisms from $\Gamma$ into $Af_{0}^{\ast}$.
Then, according to Lemma \ref{ouvR}, for small $t$,
$\rho_{t}$ is a hyperbolic representation.

Denote by ${\cal D}_{t}$ the developing maps of
the flag structures realizing the holo\-no\-my morphisms
$\rho_{t}$. They vary continuously in the compact open 
topology of maps $\widetilde{SL} \rightarrow X$.
Let $K$ be a compact fundamental domain of the action of $\Gamma$
on $\widetilde{SL}$. 
For small $t$, ${\cal D}_{t}$ is near ${\cal D}_{0}$ in the
compact--open topology, and since
${\cal D}_{0}(K)$ is a compact subset of $X_{0}$,
${\cal D}_{t}(K)$ is still a compact subset of
$X_{0}$. But the whole image of ${\cal D}_{t}$
is the $\rho_{t}(\Gamma)$--saturated of ${\cal D}_{t}(K)$,
therefore, it is contained in $X_{\infty}$.

All the discussion above shows that the deformed flag structures
we considered
are Goldman flag structures in the following meaning:

\begin{defina}
\label{defgol}
A {\em Goldman flag structure\/} is a flag structure on a 
principal Seifert ma\-ni\-fold
such that:

-\qua its holonomy morphism is a hyperbolic representation
into $Af_{0}^{\ast}$,

-\qua the image of its developing map is contained in the
open subset $X_{\infty}$.

A Goldman flag structure is {\em pure\/} if its holonomy group
does not fix a projective line.
\end{defina}

The arguments above show that Goldman flag structures
form an open subset of the space of flag structures on $M$
with holonomy group contained in $Af_{0}^{\ast}$.

We are now concerned with the problem of the existence
of Goldman flag structure on the manifold $M$ which are
not unimodular canonical flag structures.
The case of non pure Goldman flag manifolds follows from
a result of F Salein in the following way: 
consider any morphism $u$ from the cocompact
fuchsian group $\bar{\Gamma}$ into ${\bf R}^{+}$,
and consider the new subgroup $\bar{\Gamma}_{u}$ of $GL_{0}$ obtained by
replacing $\gamma$ by the multiplication of $\gamma$ by the homothety
of factor $\bar{u}(\gamma)$. The logarithm of the absolute value
of $\bar{u}$ induces a morphism $L_{u}\co  \bar{\Gamma} \rightarrow {\bf R}$.
Then:

\begin{thma}{\rm\cite{salein}}\qua
\label{sal1}
The action of $\bar{\Gamma}_{u}$ on
$X_{0}$ is free and pro\-per if and only the stable norm
of $L_{u}$ of $u$ is less than $\frac{1}{2}$.\fin
\end{thma}

\rquea
Actually, this Theorem is not stated in this form in \cite{salein}:
F Salein considered the following action of $\bar{\Gamma}$ on $PSL$: 
every element $\gamma$ maps
an element $g$ of $PSL$ on $\gamma g \Delta(\gamma)$, where $\Delta(\gamma)$
is the diagonal matrix with diagonal coefficients $e^{L_{u}(\gamma)}$,
$e^{-L_{u}(\gamma)}$. Then he proved that this action is free and proper
if and only if the stable norm of $2L_{u}$ is less than $1$ (Th\'eor\`eme $3.4$
of \cite{salein}).
But the action that we consider here is 
a double covering of the action considered by F Salein:
indeed, using the fact that $SL$ acts freely and transitively
on $X_{0}$, we identify $X_{0}$ with $SL$, and then project on $PSL$.
This double covering is an equivariant map. 
\erquea

\begin{cora}
\label{sal}
For any hyperbolic representation $\rho\co  \Gamma \rightarrow GL_{0}$,
the action of $\rho(\Gamma)$ on $X_{0}$ is free and proper.
\end{cora}

\preua
This is a corollary of Theorem \ref{sal1} and of Remark \ref{defhyp}.
Proposition \ref{propac} will give another proof of this fact.\fin

The quotients of $X_{0}$ by hyperbolic subgroups of $GL_{0}$
are called {\em canonical Goldman flag manifolds.\/}

\begin{propa}
\label{fuch}
Every hyperbolic representation is the holonomy
representation of some Goldman flag structure, which is a small
deformation of a canonical flag structure.
\end{propa}

\preua
Let $\rho\co  {\Gamma} \rightarrow Af_{0}^{\ast}$ be a hyperbolic morphism.
If $\rho(\Gamma)$ is
contained in $GL_{0}$, the proposition follows from the Corollary \ref{sal}:
the quotient of $X_{0}$ by $\rho({\Gamma})$ is a non-pure Goldman
flag manifold.

Consider now the case where $\rho(\Gamma)$ is not contained in $GL_{0}$.
Conjugating $\rho$ by a homothety of factor $s$ amounts to multiplying
the translational part of $\rho^{\ast}$ by $s$.  Therefore, if $s$
is small enough, the conjugate of $\rho$
is close to its linear part (the conjugacy does not affect
this linear part). Therefore, this conjugate is the
holonomy of some deformation of a canonical flag structure, ie,
a Goldman flag structure.
Now, conjugating back by the homothety of factor
$s^{-1}$ corresponds to multiplying the developing map
of this flag structure by $s^{-1}$. \fin

\rquea 
\label{defosa}
A corollary of Theorem B will be that
the holonomy morphism characterizes the Goldman
flag structures, ie, two Goldman flag structures
whose holonomy morphisms are conjugate in $Af_{0}^{\ast}$
are isomorphic.
As a corollary, using Proposition \ref{fuch},
Goldman flag structures are all deformations of canonical flag
structures.
\erquea

\begin{defina}
\label{golfoliation}
A Goldman foliation is the second tautological foliation
of a Goldman flag structure.
\end{defina}

\begin{propa}
\label{sa}
Goldman foliations are $Af_{0}$--foliation.
\end{propa}

\preua
As we observed previously, the second tautological foliation
of a flag manifold is transversely projective.
The holonomy morphism of this projective structure
is the holonomy morphism of the flag structure,
and its developing map is the composition of
the developing map of the flag structure with the projection $p_{2}$ of
$X$ onto ${\bf R}P^{2}$. 
For flag manifolds, the dual holonomy group is
by definition in $Af_{0}$, and the image of the developing map
is contained in $X_{\infty}$. The proposition
follows since $p_{2}(X_{\infty})$ is the affine plane ${\bf R}^{2}$.\fin

\rquea
Obvious examples of non-pure Goldman flag manifolds are
the canonical ones. They are actually the only ones.
When the holonomy group is contained in $SL$,
this follows from the proposition \ref{sa}
and from the classification of $SL$--foliations by
S Matsumoto \cite{matsu}. Theorem B provides
the proof in all the cases.
\erquea

\rquea
\label{exotic}
In the case of unimodular canonical flag structures,
the Goldman foliation is 
induced by the right action on $M$,
the left quotient of $X_{0} \approx SL$ by $\Gamma$, 
by the unipotent subgroup:
\[
\left(\begin{array}{cc}
1 & 0  \\
t & 1 
\end{array}\right)
\]
In other words, it is the horocyclic foliation of the Anosov flow
induced by diagonal matrices.

Similarly, Goldman foliations associated to non-unimo\-du\-lar
canonical Goldman flag structures are horocyclic foliations associated
to some Anosov flows: the {\em exotic Anosov flows} introduced in
\cite{ghyqf2}. Exotic Anosov flows are characterized by the following
property: they are, with the suspensions of linear hyperbolic
automorphisms of the torus, the only Anosov flows on closed
$3$--manifolds admitting a smooth splitting.  For this reason, we call
these $GL_{0}$--foliations {\em exotic horocyclic foliations\/}.
\erquea

We discuss now the problem of deformation of canonical flag
structures: what are the canonical flag structure which can be
deformed to pure Goldman flag structures?  According to the remark
\ref{defosa}, this question amounts to identifying Seifert manifolds
admitting pure Goldman structures.

The generator $h$ of the center of $\Gamma$ is mapped by
$\rho_{0}$ on the identity matrix $Id$, or its opposite
$-Id$. In the first case, $\Gamma$ is said {\em adapted,\/}
in the second one, $\Gamma$ is {\em forbidden.\/}
For example, the fundamental group of the unit tangent bundle
$M_{0}$ of $\bar{\Sigma}$ is of the forbidden type. $\Gamma$ is adapted
if and only if the finite covering $M \rightarrow M_{0}$
is of even index.

Then, $\Gamma$ admits a presentation,
with $2g+r+1$ generators $a_{i}$, $b_{i} \;\; (i=1...g)$, 
$q_{j} \;\; (j=1...r)$ and $h$,
satisfying the relations 
\[ [a_{1},b_{1}]...[a_{g},b_{g}]q_{1}...q_{r}=h^{e}, 
q_{j}^{\alpha_{j}}=h^{\beta_{j}}, [h,a_{i}]=[h,b_{i}]=[h,q_{j}]=1 
\]

\begin{propa}
\label{defor}
The canonical flag structure associated to $\Gamma$ can be
deformed to a pure Goldman flag structure if and only if $\Gamma$
is adapted.
\end{propa}

As a corollary, the canonical flag structure
on the unit tangent bundle of a hyperbolic orbifold cannot
be deformed to a pure Goldman flag structure. But its double
covering along the fibers can be deformed non-trivially.

\preuda{defor} We need to understand when the morphism $\rho_{0}$ can
be deformed in $Af^{\ast}_{0}$ to morphisms which do not preserve a
projective line.  Dually, this is equivalent to seeing when there are
morphisms $\rho\co  \Gamma \rightarrow Af_{0}$ without common fixed
point.

We first deal with the forbidden case: in this case,
the center of the holonomy group $\rho_{0}(\Gamma)$
is not trivial: it contains $-Id$. For any perturbation $\rho$,
$\rho(h)$ remains an order two element of $Af_{0}$, ie, conjugate
to $-Id$. Since $\rho(h)$ commutes with every element
of $\rho(\Gamma)$, its unique fixed point is preserved
by all $\rho(\Gamma)$. Hence, the flag structure is not pure.

Consider now the adapted case: then, $\rho_{0}(h)=Id$.
We have to find $2g$ values in $Af_{0}$ for the $\rho(a_{i})$'s and the $\rho(b_{i})$'s,
$r$ values for the $\rho(q_{j})$'s such that $\rho(q_{j})^{\alpha_{j}}=Id$,
and satisfying the relation $(\ast)$ below:
$$[\rho(a_{1}),\rho(b_{1})]...[\rho(a_{g}),\rho(b_{g})]\rho(q_{1})...\rho(q_{r})
= Id \;\;\;\;\;\;\; (\ast) $$
We realize this by adding small translation parts to the $\rho_{0}(a_{i})$,
$\rho(b_{i})$ and $\rho(q_{j})$,
ie, we try to find $\rho$ with the same linear part
than $\rho_{0}$. Adding a translationnal part to $\rho(q_{j})$
does not affect the property of being of order $\alpha_{j}$
(here $\alpha_{j}$ is bigger than $2$!)
and equation $(\ast)$ depends linearly on the
added translational parts (the linear part $\rho_{0}$ being fixed).
The number of indeterminates is $2(2g+r)$, therefore, the space of
solutions is of dimension
at least $4g+2r-2$. Amongst them, the radiants ones---ie, fixing
a point of the plane---are the
conjugates of $\rho_{0}$ by affine conjugacies whose linear parts commute
with $\rho_{0}$, ie, by compositions of homotheties and translations.
The space of radiant solutions 
is thus of dimension $3$. 
Therefore, the dimension of the space of Goldman deformations
is at least $4g+2r-5$.
But the inequality $4g+2r>5$ is always true for 
hyperbolic orbifolds. \fin

\section{Description of Goldman flag manifolds}
\label{flagol}

\begin{propa}
\label{tra}
Let $\Phi$ be an $Af_{0}$--foliation on a closed $3$--manifold $M$.
Assume that $\Phi$ is transverse to a transversely projective 
foliation $\cal F$ on $M$ of codimension one. Assume
moreover that the dual of the transverse holonomy of $\cal F$
coincides with the projectivised linear part of the
transverse holonomy of $\Phi$. Then, $\Phi$ is
the second tautological foliation of some flag structure
on $M$.
\end{propa}

\preua
Let ${\xi}\co  \widetilde{M} \rightarrow {\bf R}^{2}$ be
the developing map of the transverse structure of $\Phi$,
and $\tau\co  \widetilde{M} \rightarrow {\bf R}P^{1}$ be the
developing map for the transverse structure of $\cal F$.
Let $\rho\co  \Gamma \rightarrow Af_{0}$ be the holonomy
morphism of the transverse structure of $\Phi$.
By hypothesis, the holonomy morphism associated to $\cal F$
is the dual $\rho_{0}^{\ast}$, where $\rho_{0}$ is
the linear part of $\rho$.
Define ${\cal D}\co  \widetilde{M} \rightarrow X_{\infty}$ as
follows: for any element $m$ of $\widetilde{M}$,
${\cal D}(m)$ is the pair $(x,d)$, where $d$ is equal
to $\xi(m) \in {\bf R}^{2} \subset {\bf R}P^{2}_{\ast}$,
and where $x$ is the point of ${\bf R}P^{2}$ corresponding to
the line in ${\bf R}P^{2}$ containing $\xi(m)$ and parallel
to the direction $\tau(m)$.
Since $\Phi$ and $\cal F$ are transverse, $\cal D$
is a local homeomorphism. It is clearly
equivariant with respect to the actions on \mr and $X$ of $\Gamma$.
It is the developing map
of the required flag structure.\fin

As we will see below (section \ref{foli}), Goldman manifolds are
typical illustrations of this proposition: the $Af_{0}$--foliations
associated to a Goldman manifold are transverse to a foliation
satisfying the hypothesis of Proposition \ref{tra}.  Any small
$Af_{0}$--deformation of the $Af_{0}$--foliation (amongst the category
of $Af_{0}$--foliations), the linear part of the holonomy being
preserved, still remains transverse to the foliation. Therefore,
Proposition \ref{tra} applied in this context proves the equivalence
of the definition of Goldman foliations we have given here with the
definition introduced by Goldman, defining them as the affine
perturbations of horocyclic foliations.  Actually, the existence of
the transverse foliation is the key ingredient which allows us to
study Goldman manifolds.

Let $M$ be a Goldman flag manifold.  As usual, let $\Gamma$ be the
fundamental group of $M$, let $\cal D$ be the developing map of the
flag structure, and let $\rho\co  \Gamma \rightarrow Af_{0}^{\ast}$ be
the holonomy morphism, which is assumed to be hyperbolic.  In order to
prove Theorem B, we can replace $M$ by any finite covering of itself,
ie, replace $\Gamma$ by any finite index subgroup of itself.  In
particular, thanks to Remark \ref{notorsion}, we can assume that the
kernel of $\rho$ is $H$, and that $\rho(\Gamma)$ has no element of
finite order.

Let $\rho_{0}$ be the projectivised linear part of $\rho$. The
morphisms $\rho$ and $\rho_{0}$ induce morphisms on the surface group
$\bar{\Gamma}$, the quotient of $\Gamma$ by $H$.  We will sometimes
denote these induced morphisms abusively by $\rho$ and $\rho_{0}$.
Let $\Omega \subset X_{\infty}$ be the image of $\cal D$.

Let $\Phi$ be the Goldman foliation: it is an $Af_{0}$--foliation, 
its holonomy morphism being
$\rho$, and its developing map being ${\cal D}_{2}=p_{2} \circ {\cal D}$.
Let $\widetilde{\Phi}$ be the lifting of
$\Phi$ to the universal covering $\widetilde{M}$
of $M$.

\subsection{The affine foliation}
\label{foli}

On $X_{\infty}$, we can define the following codimension one foliation
${\cal F}_{0}$: two points $(x,d)$ and $(x',d')$ of $X_{\infty}$ are
on the same leaf if and only if there is a line containing $0$, $x$
and $x'$.  The space of leaves of ${\cal F}_{0}$ is ${\bf R}P^{1}$.
Moreover, every leaf of ${\cal F}_{0}$ is naturally equipped with an
affine structure and for this structure, the leaf is isomorphic to
the plane through the projection $p_{2}$.  The foliation ${\cal
F}_{0}$ is $Af_{0}^{\ast}$--invariant; therefore it induces a regular
foliation $\cal F$ on $M$.  Up to finite coverings, $\cal F$ is
orientable and transversely orientable.  It is a transversely
projective foliation: there is a developing map $\tau \co  \widetilde{M}
\rightarrow {\bf R}P^{1}$ and a holonomy morphism $\rho_{0}\co  \Gamma
\rightarrow PSL$.  Observe that, as our notation suggests, $\rho_{0}$
is the projectivised linear part of $\rho$.  The developing map $\tau$
is the map associating to $x$ the leaf of ${\cal F}_{0}$ containing
${\cal D}(x)$.

Let \fr denote the lifting of
$\cal F$ to $\widetilde{M}$.
Let $\cal Q$ be the leaf-space of $\widetilde{\cal F}$:
the fundamental group $\Gamma$ acts on it.

Observe that every leaf $F$ of the foliation $\cal F$ has a natural affine
structure, whose developing map is the restriction of 
${\cal D}_{2}$ to any leaf of \fr above $F$.

\begin{lema}
\label{tau}
The foliation $\cal F$ is taut, ie, ${\cal F}$ admits no Reeb component.
\end{lema}

\preua
Assume that $\cal F$ admits a Reeb component.
Let $F$ be the boundary torus of this Reeb component:
the inclusion of $\pi_{1}(F)$ in $\Gamma$ is non-injective. 
Thus, the natural affine structure
of $F$ has a non-injective holonomy morphism,
and every element of infinite order of the
holonomy group is hyperbolic.
This is in contradiction with the classification of affine
structures on the torus \cite{naganoyagi}.
\fin

It follows from a Theorem of W Thurston \cite{thuthese} that $\cal F$
is a suspension. In particular, the leaf space $\cal Q$ is
homeomorphic to the real line, and the developing map $\tau$ induces a
cyclic covering ${\cal Q} \rightarrow {\bf R}P^{1}$.  The natural
action of $\Gamma$ on the leaf space $\cal Q$ is conjugate to a
lifting of the action of the cocompact fuchsian group
$\rho_{0}(\Gamma) \subset PSL$ on ${\bf R}P^{1}$.  It follows that the
$\Gamma$--stabilizer of a point in $\cal Q$ is trivial or
cyclic. Moreover, the $\Gamma$--orbits in $\cal Q$ are dense.  In terms
of $\cal F$: every leaf of $\cal F$ is a plane or a cylinder, and is
dense in $M$.

Let $K$ be a compact fundamental domain for the action of $\Gamma$ on
$\widetilde{M}$.  Let $g$ be any $\Gamma$--invariant metric on
$\widetilde{M}$.  We fix a flat euclidian metric $dy^{2}$ on ${\bf
R}^{2}$.  This is equivalent to selecting an ellipse field $y \mapsto
\bar{E}(y)$ on the plane preserved by translations.

If the ellipses are chosen sufficiently small, the following fact is true:
for any element $x$ of $K$, there is an unique open subset $E(x)$
of the leaf through $x$ such that:

-\qua it contains $x$,

-\qua the restriction of ${\cal D}_{2} = p_{2} \circ {\cal D}$ to $E(x)$ is injective,

-\qua the image of $E(x)$ by ${\cal D}_{2}$ is $\bar{E}({\cal D}_{2}(x))$,

-\qua the $g$--diameter of $E(x)$ is less than $1$.

Since the dual morphism $\rho^{\ast}$ is hyperbolic, there exist
a real positive $\epsilon$ such that the following fact is true:
{\em for any element $\gamma$ of $\Gamma$ and for any element
$y$ of ${\bf R}^{2}$, the iterate $\rho^{\ast}(\gamma) y$ is the middle point
of an affine segment $\bar{\sigma}(\gamma y)$ of length $2\epsilon$ which is
contained in 
the ellipse $\rho^{\ast}(\gamma)\bar{E}(y)$\/} (all these metrics properties
are relative to the fixed euclidean metric $dy^{2}$).

\begin{lema}
\label{joyau}
Every leaf of $\cal F$, equipped with its affine structure,
is convex.
\end{lema}

\preua Let $\widetilde{F}$ be a leaf of $\widetilde{\cal F}$.
According to Lemma \ref{expcon}, if $\widetilde{F}$ is not convex,
there is an element $x$ of $\widetilde{F}$ for which the exponential
${\cal E}_{x}$ is not convex.  Let $U$ be the image of ${\cal E}_{x}$
by ${\cal D}_{2}$: the restriction of ${\cal D}_{2}$ to ${\cal E}_{x}$
is an affine homeomorphism over $U$. Hence, $U$ is not convex.
According to \ref{conplan}, there are two points $y$ and $z$ in ${\cal
E}_{x}$, and a closed subset $k$ of the segment $]{\cal D}_{2}(y),
{\cal D}_{2}(z)[$ such that the closed triangle $T$ with vertices
${\cal D}_{2}(x)$, ${\cal D}_{2}(y)$ and ${\cal D}_{2}(z)$ is
contained in $U$, except at $k$.  Modifying the choice of $x$ and
restricting to a smaller triangle if necessary, we can assume that the
$dy^{2}$--diameter of $T$ is as small as we want. In particular, we can
assume that for every point $y'$ sufficiently near to $k$, any segment
centered at $y'$ and of length $2\epsilon$ must intersect $[{\cal
D}_{2}(x),{\cal D}_{2}(y)] \cup [{\cal D}_{2}(x),{\cal D}_{2}(z)]]$.

Let $V$ be the subset of ${\cal E}_{x}$ that is mapped by ${\cal
D}_{2}$ to $T \setminus k$, and let $v$ be the compact subset of $V$
that is mapped onto $[{\cal D}_{2}(x),{\cal D}_{2}(y)] \cup [{\cal
D}_{2}(x),{\cal D}_{2}(z)]]Lebesgue$.  Let $\tau$ be a segment in $V$
such that ${\cal D}_{2}(\tau)$ is a segment $[{\cal D}_{2}(x), t[$,
where $t$ belongs to $k$.  Let $t_{n}$ be a sequence of points in
$\tau$ such that ${\cal D}_{2}(t_{n})$ converge to $t$. For every
index $n$, there exists an element $\gamma_{n}$ of $\Gamma$ and an
element $x_{n}$ of $K$ such that $t_{n}=\gamma_{n} x_{n}$.

We claim that the sequence $t_{n}$ escapes from any compact subset of
$\widetilde{F}$. Indeed, if this is not true, extracting a subsequence
if necessary, we can assume that $t_{n}$ converges to some point
$\bar{t}$ of $\widetilde{F}$. Clearly, ${\cal D}_{2}(\bar{t})$ is
equal to $t$.  Let $W$ be a convex neighborhood of $\bar{t}$ in
$\widetilde{F}$ such that the restriction of ${\cal D}_{2}$ to it is
injective.  According to Lemma \ref{interconv}, the restriction of
${\cal D}_{2}$ to $V \cup W$ is a homeomorphism to $T \cup {\cal
D}_{2}(W)$. It follows that the path $\tau$ can be completed as a
closed path joining $x$ to $\bar{t}$.  Hence, $\bar{t}$ belongs to
${\cal E}_{x}$, ie,  $t$ belongs to $U$. Contradiction.

Therefore the $t_{n}$ go to infinity. Their $g$--distances in
$\widetilde{F}$ to the compact set $v$ tend to infinity.  When $n$ is
sufficiently big, this distance is bigger than $1$. Therefore, none of
the ellipses $E_{n}=\gamma_{n} E(x_{n})$ intersects $v$, since their
$g$--diameter are less than $1$.  On the other hand, since ${\cal
D}_{2}(E_{n})=\rho^{\ast}(\gamma_{n})\bar{E}({\cal D}_{2}(x_{n}))$
contains the segment $\bar{\sigma}(t_{n})$ of length $2\epsilon$, the
ellipse ${\cal D}_{2}(E_{n})$ intersects $[{\cal D}_{2}(x),{\cal
D}_{2}(y)] \cup [{\cal D}_{2}(x),{\cal D}_{2}(z)]]Lebesgue$.
According to the lemma \ref{interconv}, it follows that $E_{n}$
intersects $v$. Contradiction.  \fin

In the following lemma, we call any open subset of the affine
plane bounded by two parallel lines a {\em strip\/}.

\begin{lema}
\label{jj}
The leaves of \fr are affinely isomorphic to the affine plane,
or to an affine half plane, or to a strip.
\end{lema}

\preua Let $\widetilde{F}$ be the universal covering of a leaf of
${\cal L}$.  According to Lemmas \ref{joyau} and \ref{conv}, the
restriction of ${\cal D}_{2}$ to $\widetilde{F}$ is a homeomorphism
onto a convex subset $U$ of the plane.  In order to prove the
proposition, we just have to see that the boundary components of the
convex $U$ are lines. Assume that this is not the case. Then there is
a closed half plane $P$ such that the intersection of $P$ with the
closure of $U$ is a compact convex set $K$ whose boundary is the union
of a segment $]x,y[$ contained in $U$ and a convex curve $c$ contained
in $\partial U$. We obtain a contradiction as in the proof of the
Lemma \ref{joyau} by considering ellipses centered at points $t_{n}$
of $\widetilde{F}$ such that ${\cal D}_{2}(t_{n})$ converges to some
point $t$ of $c$: for sufficiently big $n$, these ellipses, containing
segments whose length is bounded by below, must intersect $c$.  This
leads to a contradiction with the Lemma \ref{interconv}. \fin

The developing map $\tau$ induces a finite
covering of the quotient of $\cal Q$ 
by the center of $\Gamma$ over the circle ${\bf R}P^{1}$.
Let $n$ be the degree of this covering.
Consider $p_{n}\co  X_{\infty}^{n} \rightarrow X_{\infty}$, 
the finite covering of
$X_{\infty}$ of degree $n$.
Let $\widehat{M}$ be the quotient of \mr by the
center of $\Gamma$.
The map $\cal D$ induces a 
map $\widehat{\cal D}$
from $\widehat{M}$ into $X^{n}_{\infty}$.
The action of $\rho(\Gamma)$ on $X_{\infty}$ lifts
to an action of $\bar{\Gamma}$ on $X^{n}_{\infty}$ for which
$\widehat{\cal D}$ is equivariant.
Obviously, $p_{n}(\Omega^{n})=\Omega$.

\begin{propa}
\label{propre}
The map $\widehat{\cal D}$ is a homeomorphism onto
some open subset $\Omega^{n}$ of $X^{n}_{\infty}$.
\end{propa}

\preua This follows from the injectivity of ${\cal D}_{2}$ on every
leaf of \fr and from the fact that $\tau$ is a cyclic covering over
${\bf R}P^{1}$.\fin

The content of the following sections is to identify the form of
$\Omega^{n}$. It is not yet clear for example that $\Omega^{n}$ is a
cyclic covering over $\Omega$.

\subsection{Affine description of the cylindrical leaves}

Let $F_{0}$ be the lifting of a cylindrical leaf of $\cal F$.  The set
of elements of $\Gamma$ preserving $F_{0}$ is a subgroup generated by
an element $\gamma_{0}$ of infinite order.  Since $\rho(\gamma_{0})$
is a hyperbolic element of $PGL(3,{\bf R})$, $F_{0}$ is an attracting
or repelling fixed point of $\gamma_{0}$ in $\cal Q$.  We choose
$\gamma_{0}$ such that $F_{0}$ is a attracting fixed point of
$\gamma_{0}$. Observe that the fixed points of $\gamma_{0}$ in $\cal
Q$ are discrete, infinite in number, and alternatively attracting and
repelling. We denote by $F_{1}$ the lowest fixed point of $\gamma_{0}$
greater than $F_{0}$.

\begin{lema}
\label{bb}
${\cal D}_{2}(F_{0})$ is a half-plane. Its boundary $d(F_{0})$ is a line
preserved by $\rho^{\ast}(\gamma_{0})$. 
\end{lema}

\preua This follows directly from Lemma \ref{jj} and the fact that
$\rho^{\ast}(\gamma_{0})$ is a hyperbolic affine transformation acting
freely on ${\cal D}_{2}(F_{0})$. \fin

\begin{defina}
For every leaf $F$ of \fr we define 
$]F_{0},F[$ as the set of leaves in $\cal Q$ which separates $F$ from $G$.
The interval $[F_{0},F]$ is the union of $]F_{0},F[$ with
$\{ F_{0}, F \}$. We define
$\Omega_{F}$ as
the subset of elements of $F_{0}$, whose $\widetilde{\Phi}$--leaf
intersects $F$.
\end{defina}

Another equivalent definition of $[F_{0},F]$ is to consider
it as the set of leaves of \fr meeting every path
joining $F_{0}$ to $F$.

\begin{lema}
\label{omega}
The sets $\Omega_{F}$ are convex open subsets of $F_{0}$.
\end{lema}

\preua $\Omega_{F} \approx {\cal D}_{2}(\Omega_{F})$ is the
intersection of the ${\cal D}_{2}(L)$, where $L$ is in $[F_{0},F]$.
It is therefore an intersection of half-planes (maybe empty) (observe
that a strip is the intersection of two half-planes, and we can omit
the leaves $L$ whose ${\cal D}_{2}$--image are the whole plane since
they make no new contribution to the intersection).  \fin

\begin{lema}
\label{tt}
Let $F$ be an element of $]F_{0},F_{1}[$. Consider the sequence of
convex subsets of $F_{0}$, indexed by positive integers $n$, formed by
the $\gamma_{0}^{n}\Omega_{F}$.  This is an increasing sequence under
inclusion.  Moreover, the union of these convex subsets is the whole
of $F_{0}$ and the interior of their intersection is not empty.
\end{lema}

\preua Let $F'$ be a leaf such that $\Omega_{F'}$ is not empty (for
example, this is true if $F'$ is near $F_{0}$).  Since $F_{0} \prec
\gamma_{0}F' \prec F'$, then
$\gamma_{0}\Omega_{F'}=\Omega_{\gamma_{0}F'}$ contains $\Omega_{F'}$.
Since the $\gamma_{0}^{n}F'$ converge to $F_{0}$ when $n$ tend to
$+\infty$, the union of the $\gamma_{0}^{n}\Omega_{F'}$ is the whole
$F_{0}$.  Observe that for any $F$ in $]F_{0},F_{1}[$, there exists
some integer $n$ such that $\gamma_{0}^{n}F'$ is greater than $F$.
Therefore, $\Omega_{F}$ is not empty since it contains
$\gamma_{0}^{n}\Omega_{F'}$.

Since the action of $\Gamma$ on $\cal Q$ is a lifting
of the action of a cocompact fuchsian group on ${\bf R}P^{1}$, there is an element
$\gamma_{1}$ in $\Gamma$, fixing two leaves $F'_{0}$ and $F'_{1}$,
such that $]F'_{0},F'_{1}[$ contains no other fixed point
of $\gamma_{1}$, but containing $F_{0}$ and $F_{1}$.
What we did above for the pair $(\gamma_{0}, F_{0})$ can be applied to 
the pair $(\gamma_{1}, F'_{0})$:
the set of $\widetilde{\Phi}$--leaves meeting both 
$F'_{0}$ and $F_{1}$ is not empty. Since all these $\widetilde{\Phi}$--leaves
meet $F_{0}$ and $F_{1}$, $\Omega_{F_{1}}$ is not empty.
The intersection between the $\Omega_{F}$ contains $\Omega_{F_{1}}$,
Therefore, its interior is not empty.
\fin

Remember that we assumed that $F_{0}$ is an {\em attracting\/} fixed point
of $\gamma_{0}$.

\begin{cora}
\label{boundary}
The boundary line of $F_{0}$ is the unstable
line of $\rho^{\ast}(\gamma_{0})$, ie, it is parallel
to the eigenspace associated to the eigenvalue of
$\rho^{\ast}(\gamma_{0})$ of absolute value greater than $1$.
\end{cora}

\preuda{boundary} Assume that the Lemma is false.  Take some leaf $F$
in $]F_{0},F_{1}[$.  According to Lemma \ref{tt}, the
$\gamma_{0}^{n}\Omega_{F}$ for positive $n$ form an increasing
sequence of convex sets whose union is the whole of $F_{0}$.  This is
possible only if the convex set $\Omega_{F}$ is a strip containing
$d(F_{0})$ in its boundary. But then the intersection of the
$\gamma_{0}^{n}\Omega_{F}$ would be empty: this contradicts
\ref{tt}.\fin

\begin{cora}
\label{strip}
No leaf of \fr is a strip.
\end{cora}

\preua Assume that some leaf $F$ is a strip. Then it admits at least
two iterates $\gamma F$ and $\gamma' F$ in $]F_{0}, F_{1}[$.  One of
them, let's say $\gamma F$, disconnects $F_{0}$ from the other
($\gamma' F$).  We can choose these iterates such that the strips
${\cal D}_{2}(\gamma F)$ and ${\cal D}_{2}(\gamma' F)$ are not
parallel.  Then, the intersection of these two strips is a
parallelogram.  But, since $\gamma F$ disconnects $F_{0}$ from $\gamma'
F$, this parallelogram contains ${\cal D}_{2}(\Omega_{\gamma' F})$,
According to Lemma \ref{tt}, the intersection between the positive
$\rho^{\ast}(\gamma_{0})$--iterates of this parallelogram must have a
non-empty interior.  But this is clearly impossible: for any
parallelogram $P$ of the plane, the intersection of the positive
$\rho^{\ast}(\gamma_{0})$--iterates of $P$ is either empty, either a
subinterval of the unstable line of $\rho^{\ast}(\gamma_{0})$. \fin

\subsection{Description of the image and the limit set}
\label{calU}

Let $U(\rho)$ be the image of ${\cal D}_{1}=p_{1} \circ {\cal
D}$. This is a subset of ${\bf R}P^{2} \setminus 0$.  Let ${\cal
G}_{0}$ be the projection in ${\bf R}P^{2}$ of ${\cal F}_{0}$: this is
the foliation whose regular leaves are the $d \setminus 0$, where $d$
is any projective line in ${\bf R}P^{2}$ containing $0$.

\begin{lema}
\label{pivot}
$U(\rho)$ is not the whole of ${\bf R}P^{2} \setminus 0$.
\end{lema}

\preua Let $\gamma_{0}$ be an element of $\Gamma$ of infinite order
admitting fixed points in $\cal Q$.  Let $d_{0}$ be the attracting
fixed point of $\rho(\gamma_{0})$ in ${\bf R}P^{2}$, and $g_{0}$ the
leaf of ${\cal G}_{0}$ containing $d_{0}$: in ${\bf R}P^{1}$, the
space of leaves of ${\cal G}_{0}$, $g_{0}$ is an attracting fixed
point of $\rho_{0}(\gamma_{0})$. Since the action of $\gamma_{0}$ on
$\cal Q$ is a lifting of the action of $\rho_{0}(\gamma_{0})$ on ${\bf
R}P^{1}$ and since this action admits a fixed point, every leaf of
$\widetilde{\cal F}$ whose image is contained in $g_{0}$ is an
attracting fixed point of $\gamma_{0}$.  According to Corollary
\ref{boundary}, the image by ${\cal D}_{2}$ of such a leaf is an
half-plane whose boundary is the unstable line of
$\rho^{\ast}(\gamma_{0})$, ie $d_{0}$. It follows that $d_{0}$ is
not in $U(\rho)$.  \fin

There is a morphism $\rho_{1}\co  \Gamma \rightarrow GL_{0}$ and
two maps $u$ and $v$ from $\Gamma$ into ${\bf R}$ such that
the morphism $\rho\co  \Gamma \rightarrow PGL(3,{\bf R})$ is induced 
by a morphism of the form:
\[
\rho(\gamma) = \left(\begin{array}{cc}
\rho_{1}(\gamma) & \begin{array}{c}
     0 \\
     0
    \end{array}\\
\begin{array}{cc}
u(\gamma) & v(\gamma) 
\end{array} & 1
\end{array}\right)
\]
The open set $U(\rho)$ is a set of lines in ${\bf R}^{3}$.  Their
union minus the origin is $\rho(\Gamma)$--invariant open cone in ${\bf
R}^{3}$. We denote this by ${\cal U}(\rho)$.

For any point $w=(x,y)$ in ${\bf R}^{2} \setminus 0$, consider
the set of real numbers $z$ such that $(x,y,z)$ belongs to ${\cal U}(\rho)$.
We denote it by $I(w)$; it is an open subset of ${\bf R}$.
According to Lemma \ref{jj} and Corollary \ref{strip}, for every leaf $F$
of \fr, the image ${\cal D}_{2}(F)$ 
is the whole plane or a half-plane. In both cases, $d_{\infty}$
meets the boundary of ${\cal D}_{2}(F)$. It means that
the point $0$ is an extremity of ${\cal D}_{1}(F)$. It follows that
that for every $w$ in ${\bf R}^{2} \setminus 0$, 
$I(w)$ contains an interval of the form $]-\infty, t[$, and
another of the form $]t, +\infty[$. We denote by $\delta^{-}(w)$
(respectively $\delta^{+}(w)$) the supremun (respectively the infimum) of
the real numbers $t$ for which $]-\infty, t[$ (respectively $]t,+\infty[$)
is contained in $I(w)$. These maps have the following properties:

-\qua $\delta^{-}(w) \leq \delta^{+}(w)$ or $\delta^{-}(w)=+\infty=-\delta^{+}(w)$,

-\qua $\delta^{+}(w)=-\delta^{-}(-w)$ (because $-{\cal U}(\rho)={\cal U}(\rho)$),

-\qua $\delta^{-}$ is lower semi-continuous (l.s.c.), and $\delta^{+}$ is 
upper semi-continuous (u.s.c.)
(because ${\cal U}(\rho)$ is open),

-\qua $\delta^{+}$ and $\delta^{-}$ are homogeneous of degree $1$ (because
${\cal U}(\rho)$ is a cone).

-\qua since ${\cal U}(\rho)$ is $\rho(\Gamma)$--invariant,
for every element $\gamma$ of $\Gamma$, they both satisfy:
\[ \delta^{\pm}(\rho_{1}(\gamma)(x,y))=
\delta^{\pm}(x,y) + u(\gamma)x+v(\gamma)y \;\;\;\; (\ast)\]

\begin{propa}
\label{delta}
The maps $\delta^{+}$ and $\delta^{-}$ are equal and take only finite
values.
\end{propa}

\preua For any $w$ in ${\bf R}^{2} \setminus 0$, let $\Delta(w)$ be
the difference $\delta^{+}(w)-\delta^{-}(w)$. The map $\Delta\co  {\bf
R}^{2} \setminus 0 \rightarrow {\bf R} \cup \{-\infty\}$ is u.s.c.\
and, according to equation $(\ast)$ above, $\Delta$ {\em is
$\rho_{1}(\Gamma)$--invariant.\/} On the other hand, the quotient of
$X_{0}$ by $\rho_{1}(\Gamma)$ is a canonical flag manifold
$M(\rho_{1})$, whose second tautological foliation is an exotic
horocyclic foliations (see Remark \ref{exotic}).  Therefore, $\Delta$
induces an u.s.c.\ function $\bar{\Delta}\co  M(\rho_{1}) \rightarrow
{\bf R} \cup \{ -\infty \}$ which is invariant along the leaves of the
exotic horocyclic foliation.  Since $M(\rho_{1})$ is compact, and
since $\bar{\Delta}$ is u.s.c.\ $\bar{\Delta}$ attains its maximal
value and the locus where it attains this maximal value is closed.
Since horocyclic foliations of Anosov flows are notoriously minimal
(when the flow is not a suspension, and this is the case here; see
Theorem 1.8 of \cite{plante3}), it follows that $\Delta$ is
constant. Since it is homogeneous of degree $1$, the constant value of
$\Delta$ is either $-\infty$ or $0$.  If the constant value is
$-\infty$, then $\delta^{+}$ and $\delta^{-}$ have infinite value
everywhere. Then, all the $(x,y,z)$, for $(x,y)$ describing ${\bf
R}^{2} \setminus 0$, belong to ${\cal U}(\rho)$. This contradicts
Lemma \ref{pivot} and proves the proposition.  \fin

\begin{cora}
\label{conjug}
The actions of $\rho(\Gamma)$ and $\rho_{1}(\Gamma)$ on
the projective plane are topologically conjugate, ie, there is
a homeomorphism $f\co  {\bf R}P^{2} \rightarrow {\bf R}P^{2}$
satisfying the following equivariance property:
\[ \forall \gamma \in \Gamma \;\;\; 
f \circ \rho_{1}(\gamma) = \rho(\gamma) \circ f\]
Moreover, the conjugacy $f$ is unique up to composition on
the left by homotheties.
\end{cora}

\preua
According to the Proposition \ref{delta}, the map $\delta^{+}=\delta^{-}$ is
continuous, since it is u.s.c.\ and l.s.c.\ at the same time.
Consider the following map of
${\bf R}^{3}$ minus the $z$--axis into itself:
\[ (x,y,z) \mapsto (x,y,z+\delta^{+}(x,y))\]
Since $\delta^{+}$ is homogeneous of degree one, and since
$\delta^{+}(-w)=-\delta^{-}(w)=-\delta^{+}(w)$, this map induces a
homeomorphism of ${\bf R}P^{2} \setminus \{ 0 \}$ onto itself. This
homeomorphism extends as a homeomorphism $f$ of ${\bf R}P^{2}$ onto
itself by setting $f(0)=0$.  Equation $(\ast)$ above implies the
required $\Gamma$--equivariance of $f$.  If $f'$ is another topological
conjugacy, then $f'^{-1} \circ f$ is a transformation of ${\bf R}^{2}
\setminus 0$ commuting with the linear action of
$\rho_{1}(\Gamma)$. Then, according to the rigidity of horocyclic
flows, $f'^{-1} \circ f$ must be a homothety (see \cite{abe} for the
case of geodesic flows, the case of exotic Anosov flows is similar).
\fin

Let $\Lambda(\rho)$ be the image by $f$ of the $GL_{0}$--invariant
projective line.  This is a Jordan curve.

\begin{cora}
\label{Uh}
The curve $\Lambda(\rho)$ is the closure of the union of the repelling fixed
points of elements of $\rho(\Gamma)$. 
It is the complement in
${\bf R}P^{2}$ of the disc $U(\rho) \cup \{ 0 \}$. The action
of $\rho(\Gamma)$ on $\Lambda(\rho)$ is topologically conjugate
to the projective action of the fuchsian group $\lambda(\rho)(\Gamma)$
on the projective plane ${\bf R}P^{1}$.
The action
of $\rho(\Gamma)$ on $U(\rho)$ is uniquely ergodic.
\end{cora}

\preua Using the equivariant map $f$, it is enough to check all these
statements in the case of canonical flag manifolds, in which they are
easily established.  \fin

\begin{lema}
\label{linel}
Two hyperbolic actions of $\Gamma$ are topologically conjugate if and only
if their linear parts are conjugate in $GL_{0}$.
The curve $\Lambda(\rho)$ is a projective line if and only if
the conjugacy $f$ with the linear part is projective.
\end{lema}

\preua The first part is a corollary of the rigidity of exotic
horocyclic flows. For the second part, when $\Lambda(\rho)$ is a
projective line, there is a projective transformation $g$ mapping
$U(\rho)$ to ${\bf R}^{2} \setminus {0}$, and thus mapping the
$\rho$--action of $\Gamma$ to some linear action. Then, $g \circ f$ is
a topological conjugacy between two linear actions. By the first part,
by modifying $g$, we can assume that $g \circ f$ commutes with the
linear action of $\rho_{1}(\Gamma)$ on ${\bf R}^{2} \setminus 0$.
By Lemma \ref{conjug} it is a homothety; therefore, $f$ is
projective.\fin

\begin{lema}
\label{diff}
If the map $\delta^{+}=\delta^{-}$ is differentiable on a set of non
zero Lebesgue measure, then the conjugacy $f$ is a projective
transformation.
\end{lema}

\preua The idea of the proof is due to A Zeghib.  In the hypothesis of
the lemma, since the action of $\rho_{1}(\Gamma)$ on ${\bf R}R^{2}
\setminus 0$ is uniquely ergodic, $\delta^{+}$ is differentiable
almost everywhere. We can then define an equivariant measurable map
$\tau\co  {\bf R}P^{1} \rightarrow {\bf R}^{2} \subset {\bf
R}P^{2}_{\ast}$ defined almost everywhere, by associating to every
$[x;y]$ the projective line tangent to $\Lambda(\rho)$ at the ray
$\bar{f}([x;y])=[x;y;\delta^{+}(x,y)]$: observe that $\tau([x;y])$
never contains $0$.  Let $\cal P$ be the product ${\bf R}P^{1} \times
{\bf R}P^{1}$ minus the diagonal.  Observe that the diagonal action of
$\rho(\Gamma)$ on $\cal P$ admits an ergodic invariant measure
equivalent to the Lebesgue measure.  We say that a subset of $\cal P$
is conull if the measure of its complement in $\cal P$ is $0$.  The
crucial and classical observation is that this ergodicity property
implies that there is no measurable equivariant map from $\cal P$ into
a topological space where $\Gamma$ acts freely and properly
discontinuously.

Assume that the set of pairs $(\theta,\theta')$ for which $\theta'$
does not belong to $\tau(\theta)$ is conull. Then, its intersection
with its image by the flip map $(\theta,\theta') \mapsto
(\theta',\theta)$ is conull, and its intersection with all its
$\bar{\Gamma}$--iterates also.  Thus, there is a conull
$\rho(\Gamma)$--invariant subset $\cal E$ of $\cal P$ of pairs
$(\theta, \theta')$ for which the projective lines $\tau(\theta)$ and
$\tau(\theta')$ intersect at some point $x(\theta,\theta')$ different
from $\bar{f}(\theta)$ and $\bar{f}(\theta')$.  We have then two
cases: either almost every $x(\theta,\theta')$ belongs to
$\Lambda(\rho)$ or almost all of them belongs to $U(\rho)$.  In the
first case, we obtain a $\rho(\Gamma)$--equivariant map from $\cal E$
into the set of distinct triples of points of ${\bf R}P^{1}$.  Since
the action of $\bar{\Gamma}$ on this set of triples is free and
properly discontinuous, we obtain a contradiction with the ergodic
argument discussed above. In the second case, the map associating to a
pair $(\theta,\theta')$ the flag $(x(\theta,\theta'),\tau(\theta))$ is
an equivariant map from $\cal E$ into $X(\rho)$. According to
Proposition \ref{propac} below, we obtain once more a contradiction
with the ergodic argument.

Therefore, the measure of the set of pairs $(\theta,\theta')$
for which the line $\tau(\theta)$ contains $\theta'$ is conull. 
Then, by Fubini's Theorem, there is an element $\theta$ of ${\bf R}P^{1}$
such that for almost all $\theta'$ in ${\bf R}P^{1}$, $\bar{f}(\theta')$
belongs to $\tau(\theta)$. But the intersection of $\Lambda(\rho)$ with
$\tau(\theta)$ is closed, and $\bar{f}$ is continuous: it follows
that $\Lambda(\rho)$ must be equal to $\tau(\theta)$. 
We conclude by applying Lemma \ref{linel}.
\fin

\begin{cora}
\label{lip}
The Jordan curve $\Lambda(\rho)$ is Lipschitz if and only if it
is a projective line, ie, if and only if the conjugacy $f$ is projective.
\end{cora}

\preua this follows from Lemma \ref{diff} since Lipschitz maps are
differentiable almost everywhere. \fin

\subsection{Properness of the action}

We define 
$X(\rho)$ as the intersection of $X_{\infty}$
with the preimage by $p_{1}$ of $U(\rho)$.

\begin{propa}
\label{propac}
The action of $\rho(\Gamma)$ on $X(\rho)$ is free and properly
discontinuous.
\end{propa}

\preua
The action of $\rho(\Gamma)$ on $U(\rho)$ is conjugate to the
action of $\rho_{1}(\Gamma)$ on the punctured affine plane.
Therefore, it is free, and the action of $\rho(\Gamma)$ on
$X(\rho)$ is free. 
Remember also that by replacing $\Gamma$ by a finite index subgroup, we can assume
that all the non-trivial elements of $\rho(\Gamma)$ are hyperbolics.

Since $\rho(\Gamma)$ is discrete in $Af_{0}^{\ast}$, 
we just have to establish the properness of its action on $X(\rho)$.
Assume {\em a contrario\/} that it is not the case:
there are elements $(x_{n},d_{n})$, $(x'_{n},d'_{n})$
of $X(\rho)$, and elements $g_{n}$ of $\rho(\Gamma)$
such that:

-\qua $(x'_{n},d'_{n}) = g_{n}(x_{n},d_{n})$, 

-\qua the $(x_{n},d_{n})$ converge to some $(x,d)$ in $X(\rho)$,

-\qua the $(x'_{n},d'_{n})$ converge to some $(x',d')$ in $X(\rho)$,

-\qua the $g_{n}$ escape from any compact subset of $Af_{0}^{\ast}$.

Define $g^{\ast}_{n}=\theta(g_{n})$: they escape from any
compact subset of $Af_{0}$ too.
As elements of $Af_{0} \subset PGL(3,{\bf R})$, the $g_{n}^{\ast}$ are representated
by $3 \times 3$--matrices of the
form:
\[
\left(\begin{array}{cc}
B_{n} & \begin{array}{c}
      u_{n} \\
      v_{n}
    \end{array} \\
\begin{array}{cc} 0 & 0
\end{array} & 1
\end{array}\right)
\]
For any vector subspace $E$ of ${\bf R}^{3}$ (or its dual),
we denote by $S(E)$ its projection in 
${\bf R}P^{2}$ (or ${\bf R}P^{2}_{\ast}$).
We see $GL(3,{\bf R})$ as a subset of $M(3,{\bf R})$, the algebra
of $3 \times 3$--matrices. Denote by $\Vert_{0}$ the operator norm
on $M(3,{\bf R})$; let $\cal B$ be the unit ball of
this norm. 

Extracting a subsequence if necessary, we can assume that
the sequences $h^{n}=\frac{g_{n}}{\Vert g_{n}\Vert_{0}}$ and
$h^{\ast}_{n}=\frac{g^{\ast}_{n}}{\Vert g^{\ast}_{n}\Vert_{0}}$ 
converge respectively 
to $\bar{g}$ and $\bar{g}^{\ast}$ in $\cal B$. 

A fundamental fact is the following claim: {\em the norm $\Vert g_{n}^{\ast} \Vert_{0}$
tends to $+\infty$.\/} Indeed: remember the discussion in 
Remark \ref{defhyp}. The linear part $B_{n}$ is of the form 
$\bar{u}(\gamma_{n})\rho_{0}^{\ast}(\gamma_{n})$, where $\rho_{0}^{\ast}\co  \Gamma \rightarrow SL$
is the composition of the linear part of $\rho^{\ast}\co  \Gamma \rightarrow Af^{\ast}_{0}$
with the projection of $GL_{0}$ over $SL$, and
$\bar{u}\co  \Gamma \rightarrow {\bf R}^{+}$ is a morphism.
The projection of $\rho_{0}^{\ast}(\gamma)$ in $PSL$ is the projectivised
linear part $\lambda(\rho^{\ast})(\gamma)$. Let $t(\gamma)$ be the logarithm
of the spectral radius of $\rho_{0}^{\ast}(\gamma)$: it
is also the logarithm of the spectral radius of $\lambda(\rho^{\ast})(\gamma)$. 
Since $\lambda(\rho^{\ast})(\Gamma)$
is a cocompact fuchsian group, $t(\gamma_{n})$ tends to $+\infty$ when
$n$ goes to infinity. Let now $L_{u}(\gamma)$ be the logarithm of the absolute
value of $\bar{u}(\gamma)$. Since $\rho^{\ast}$ is hyperbolic, the
stable norm of the morphism induced on $\bar{\Gamma}$ by $L_{u}$ is
less than $\frac{1}{2}$; let $0 < C < 1$ be the double of this norm: 
by definition of the stable norm,
the absolute value of $L_{u}(\gamma_{n})$ is less than $Ct(\gamma_{n})$.
It follows that $ t(\gamma_{n}) \pm L_{u}(\gamma_{n})$ is
bigger than $(1-C)t(\gamma_{n})$, and thus, that 
$t(\gamma_{n}) \pm L_{u}(\gamma_{n})$ tends to $+\infty$ with $n$.
But the absolute value of the eigenvalues of $B_{n}$ are the
exponentials of $t(\gamma_{n}) \pm L_{u}(\gamma_{n})$ and
of $-t(\gamma_{n}) \pm L_{u}(\gamma_{n})$. It follows that
one of these eigenvalues tends to $+\infty$, and therefore, that
the norm of $B_{n}$ tends to $+\infty$.

Hence, $\bar{g}^{\ast}$ is of the form:
\[
\left(\begin{array}{cc}
B & \begin{array}{c}
      u \\
      v
    \end{array} \\
\begin{array}{cc} 0 & 0
\end{array} & 0
\end{array}\right)
\]
Therefore, 
the image $I^{\ast}$ and the kernel $K^{\ast}$ of $\bar{g}^{\ast}$ are proper 
subspaces of ${\bf R}^{3}$,
and $S(I^{\ast})$ is contained in the
line at infinity.

Similar considerations show that the norm of $g_{n}$ tends to $+\infty$,
that the image $I$ and the kernel $K$
of $\bar{g}$ are proper subspaces, and that 
$S(K)$ contains the point $0$.

For every index $n$, the products $h_{n}^{t}h_{n}^{\ast}$ and
$h_{n}^{\ast}h_{n}^{t}$ (where $h_{n}^{t}$ is the transposed matrix of
$h_{n}$) are both equals to $\frac{id}{\Vert g_{n}\Vert_{0} \Vert
g_{n}^{\ast}\Vert_{0}}$.  Hence, when $n$ goes to infinity, we obtain:
\[ \bar{g}^{t} \circ \bar{g}^{\ast} = 0 = \bar{g}^{\ast} \circ \bar{g}^{t} \]
The transposed matrix $\bar{g}^{t}$ of $\bar{g}$ has to be considered
as a linear endomorphism of the dual of ${\bf R}^{3}$: $\bar{g}^{t}$
maps a linear form $\varphi$ on $\varphi \circ \bar{g}$.  The elements
of its kernel are the linear forms whose kernels contain $I$, and the
elements of its image are the linear forms whose kernels contain $K$.
From the equalities above, we obtain that {\em any element $d$ of
$S(I^{\ast})$, viewed as a projective line in ${\bf R}P^{2}$, contains
$S(I)$, and that any element $d$ of ${\bf R}P^{2}_{\ast}$ containing
$S(K)$ necessarily belongs to $S(K^{\ast})$.\/}

Observe that $\bar{g}$ (respectively $\bar{g}^{\ast}$) defines a map
from ${\bf R}P^{2} \setminus S(K)$ (respectively ${\bf R}P^{2}_{\ast}
\setminus S(K^{\ast})$) into $S(I) \subset {\bf R}P^{2}$ (respectively
$S(I^{\ast}) \subset {\bf R}P^{2}_{\ast})$.  We claim: {\em $g_{n}$
converges on ${\bf R}P^{2} \setminus S(K)$ to the map $\bar{g}$. This
convergence is uniform on the compact subsets of ${\bf R}P^{2}
\setminus S(K)$.\/} Indeed: consider any compact subset $\cal K$ in
${\bf R}P^{2} \setminus S(K)$.  We denote by $\Vert$ the euclidean
norm of ${\bf R}^{3}$.  Observe that $g_{n}$ and
$h_{n}=\frac{g_{n}}{\Vert g_{n}\Vert_{0}}$ have the same actions on
${\bf R}P^{2}$.  For every element $x$ of ${\bf R}^{3}$ representing
an element of $\cal K$ we have:
\[
\begin{array}{rcl}
\Vert \frac{h_{n}(x)}{\Vert h_{n}(x)\Vert} - \frac{\bar{g}(x)}{\Vert \bar{g}(x)\Vert}\Vert & \leq &
\Vert \frac{h_{n}(x)}{\Vert h_{n}(x)\Vert}-\frac{h_{n}(x)}{\Vert \bar{g}(x)\Vert }\Vert
+ \Vert \frac{h_{n}(x)}{\Vert \bar{g}(x)\Vert } - \frac{\bar{g}(x)}{\Vert \bar{g}(x)\Vert }\Vert \\\phantom{\vrule height .55cm depth .3cm}
  & \leq  &
\frac{\mid \Vert h_{n}(x)\Vert -\Vert \bar{g}(x) \Vert \mid }{\Vert \bar{g}(x)\Vert} +
 \frac{\Vert h_{n}(x) - \bar{g}(x)\Vert}{\Vert \bar{g}(x)\Vert} \\
 & \leq &
\frac{2}{\Vert \bar{g}(x)\Vert}\Vert h_{n}(x)-\bar{g}(x)\Vert
\end{array}
\]
The claim now follows from the fact that $\Vert \bar{g}(x) \Vert$ is
bounded from below by a positive constant valid for all the points
of the unit ball of ${\bf R}^{3}$
representing elements of $\cal K$.

The similar
property for ${g}_{n}^{\ast}$ is also true.

It follows that if $d$ does not belong to $S(K^{\ast})$, then $d'$
belongs to $S(I^{\ast})$. But this is impossible since $S(I^{\ast})$
is contained in the line at infinity.  Hence, $d$ belongs to
$S(K^{\ast})$.  Assume that $S(K^{\ast})$ is reduced to $\{ d
\}$. Then, according to the property of uniform convergence, we see
that for any small disc $D$ in ${\bf R}P^{2}_{\ast}$ containing $d$,
the closure of $D$ is contained in $g_{n}^{\ast}D$ when $n$ is
sufficiently big. It follows that $g_{n}^{\ast}$ admits a repelling
fixed point in $D$: contradiction.

Thus, $S(K^{\ast})$ is a projective line, and $S(I^{\ast})$ a single
point.  Consider now $S(K)$: it contains $0$.  We have seen that any
projective line containing $S(K)$ belongs to $S(K^{\ast})$: therefore,
if $S(K)$ was reduced to $\{ 0 \}$, $S(K^{\ast})$ would be the line
$d_{\infty}$: this is impossible since $d$ belongs to $S(K^{\ast})$ and
not to $d_{\infty}$.  Therefore, $S(K)$ is a projective line, and
$S(I)$ a single point.  The Jordan curve $\Lambda(\rho)$ does not
contain $0$, therefore, it is not $S(K)$.  It follows that there is a
point $x_{0}$ of $\Lambda(\rho)$ which does not belong to $S(K)$. The
iterates $g_{n}x_{0}$ all belong to $\Lambda(\rho)$ and converge to
$S(I)$: it follows that $S(I)$ belongs to $\Lambda(\rho)$.  Hence,
$x'$ is not $S(I)$, this implies that $x$ belongs to $S(K)$.  We use
once more the fact that any projective line containing $S(K)$ belongs
to $S(K^{\ast})$: this shows that $S(K)$, as a point in ${\bf
R}P^{2}_{\ast}$, belongs to $S(K^{\ast})$.  Dually, $S(K^{\ast})$, as
a point in ${\bf R}P^{2}$, belongs to $S(K)$. Then, $d$ and $S(K)$
have two distinct points in common: $x$ and the point
$S(K^{\ast})$. They must be equal.  But $S(K)$ contains $0$, and, by
hypothesis, $d$ does not contain $0$. This is a contradiction.  \fin

Since $U(\rho)$ is topologically an annulus, $X(\rho)$ is homeomorphic
to ${\bf S}^{1} \times {\bf R}^{2}$.  It follows that the quotient
$M(\rho)$ of $X(\rho)$ by $\rho(\Gamma)$ is a $K(\Pi,1)$. We deduce
from homological considerations that $M(\rho)$ is a {\em compact\/}
$3$--manifold.  Now, $\cal D$ induces a local homeomorphism of $M$ in
$M(\rho)$.  Since both are compact $3$--manifolds, this induced map is
a finite covering. We have proved Theorem B.

\section{The tautological foliations}
\label{flots}

In this section, we study the tautological foliations associated to
Goldman flag structures.  We are only interested in dynamical
properties which are not perturbed by finite coverings. Therefore, we
can, and we do, assume that $M$ is the quotient of $X(\rho)$ by
$\rho(\Gamma)$, where $\rho$ is the holonomy morphism.  Therefore, the
holonomy group is isomorphic to $\bar{\Gamma}$, the quotient of the
fundamental group $\Gamma$ by its center $H$. From now on, we denote
$\rho(\Gamma)$ by $\bar{\Gamma}$. We can assume that $\bar{\Gamma}$
has no torsion.

We call $\Psi$ the first tautological foliation, and $\Phi$ the second
one.  We will see that their dynamical behaviors are quite different.
We call their liftings in the covering $X(\rho)$ of $M$, $\hat{\Psi}$
and $\hat{\Phi}$.  Observe that these foliations are orientable since
$Af_{0}$ preserves any orientation of ${\bf R}^{2}$.

\subsection{Study of the first tautological foliation}

We need to consider another foliation on $M$: as a topological
manifold, $M$ is homeomorphic to the left quotient of $SL$ by the
linear part $\bar{\Gamma}_{0} \subset GL_{0}$ of $\bar{\Gamma}$.  On
this quotient, which we denote by $M_{l}$, we have the horocyclic flow
$\Psi_{0}^{t}$, induced by the right action of unipotent matrices.
Observe that the operation of ``taking the linear part" defines an
isomorphism $\bar{\Gamma} \rightarrow \bar{\Gamma}_{0}$.

\begin{thma}
\label{psi}
${\Psi}$ is topologically conjugate to
the horocyclic foliation ${\Psi}_{0}$.
\end{thma}

\preua According to Corollary \ref{Uh}, there is a topological
conjugacy $f$ between the action of $\bar{\Gamma}_{0}$ on ${\bf R}^{2}
\setminus \{ 0 \}$ and the action of $\bar{\Gamma}$ on $U(\rho)$.  But
the pairs $({\bf R}^{2} \setminus \{ 0 \}, \bar{\Gamma}_{0})$ and
$(U(\rho),\bar{\Gamma})$ can be interpreted as the leaf spaces of
$\Psi_{0}$ and $\Psi$ respectively.  The holonomy covering of the
leaves of $\Psi$ and $\Psi_{0}$ are contractible. Hence,
$(M_{l},\Psi_{0})$ and $(M,\Psi)$ are representatives of the
classifying spaces of their transverse holonomy groupoids.  $({\bf
R}^{2} \setminus \{ 0 \}, \bar{\Gamma})$ and $(U(\rho),\bar{\Gamma})$
are also representatives of these classifying spaces. Since they are
conjugate, and by uniqueness of classifying spaces modulo
equivalence, {\em there exists a homotopy equivalence $F\co  M_{l}
\rightarrow M$ mapping every leaf of $\Psi_{0}$ into a leaf of
$\Psi$.\/} Moreover, $F$ lifts to some mapping $\hat{F}$ between the
coverings $X_{0}$ and $X(\rho)$, which induces $f$ at the level of the
leaf spaces. In particular, $F$ maps two different leaves of
$\Psi_{0}$ into two different leaves of $\Psi$ (for the notion of
classifying spaces, and for all the arguments used here, we refer to
\cite{haefliger}). The problem is that this map has no reason to be
injective along the leaves of $\Psi_{0}$.

We will modify $F$ along the leaves of $\Psi_{0}$ in order to correct
this imperfection.  This idea of diffusion process along the leaves
seems due to M Gromov.  It has been used in \cite{bar1}, \cite{matsu},
and previously in \cite{ghys}.

First, we choose arbitrary parametrisations
$\Psi_{0}^{t}$ and $\Psi^{s}$ of the foliations. 
Since $\Psi_{0}^{t}$ has no periodic orbit, we
have a continuous map $u\co  M_{l} \times {\bf R} \rightarrow {\bf R}$
satisfying:
\[ \forall t \in {\bf R} \;\;\; \forall x \in M \;\;\;\;\; 
F(\Psi_{0}^{t} x) = \Psi^{u(t,x)}(F(x)) \]
$u$ is a cocycle, ie, for every element $x$ of $M_{l}$:
\[ \begin{array}{crcl}
\forall s,t \in {\bf R} \;\; & \; u(t+s,x) & = & u(t,\Psi_{0}^{s}x)+u(s,x) \\
\forall t \in {\bf R} \;\; & \; u(0,x) & = & 0
\end{array} \]
The main lemma is:

\begin{lema}
\label{existt}
There is a real $T>0$ such that, for any element $x$ of $M_{l}$, 
the quantity $u(T,x)\/$ is not zero.
\end{lema}

Assume that Lemma \ref{existt} is true.  Let $T$ be the real given by
the lemma. We define:
\[ \forall x \in M_{l} \;\;\;\; 
u_{T}(x) = \frac{1}{T} \int_{0}^{T}{u(s,x)\,ds} \]
Then, we define $F_{T}\co M_{l} \rightarrow M$:
\[ F_{T}(x)=\Psi^{u_{T}(x)}(F(x))  \]
This map has the same properties than $F$. Moreover
\[ F_{T}(h^{t}(x))=\Psi^{v_{T}(t,x)}(F_{T}(x)) \]
where:
\[ v_{T}(t,x)=\frac{1}{T} \int_{t}^{T+t}{u(s,x)\,ds} \]
The derivation of $v_{T}$ with respect to $t$ is:
\[ \begin{array}{rcl}
 \frac{\partial}{\partial t}v_{T}(t,x) & = & \frac{1}{T} [u(T+t,x) -u(t,x)] \\
 \phantom{\vrule height .5cm} & = & \frac{1}{T} u(T,\Psi_{0}^{t}(x))
\end{array} \]
According to our choice of $T$, this is never zero.  It follows that
$F_{T}$ is injective along the leaves of $\Psi_{0}$, and therefore
injective. Since it is a homotopy equivalence, it is a topological
conjugacy between $\Psi_{0}$ and $\Psi$.  Therefore, in order to
finish the proof of Theorem \ref{psi}, we just have to prove
\ref{existt}:

\preuda{existt}
Assume that Lemma \ref{existt} is not true.
Then, there is a sequence of increasing real numbers
$t_{n}$, converging to $+\infty$, and
a sequence of points $x_{n}$ in $M_{l}$ such that
the $u(t_{n},x_{n})$ are zero.
We can assume that $x_{n}$ converges to some point of $M_{l}$.
Remember that $SL$ is naturally identified with $X_{0}$.
The $x_{n}$ lift in $X_{0}$ to pairs $(y_{n},d_{n})$
where the $y_{n}$ are points in ${\bf R}^{2} \setminus \{ 0 \}$
and the $d_{n}$ are projective lines through $y_{n}$ (but not $0$) converging
to some element $(\bar{y}, \bar{d})$ of $X_{0}$.
The $\Psi_{0}^{t_{n}}(x_{n})$ lift to pairs $(y_{n},d'_{n})$.
Since $t_{n}$ go towards infinity, and since the $y_{n}$ converge
to $\bar{y}$, the $d'_{n}$ converge to the projective line containing
both $0$ and $\bar{y}$.
Now, the nullity of $u(x_{n},t_{n})$ means that
the $\hat{F}(y_{n},d_{n})$ and $\hat{F}(y_{n},d'_{n})$
are equal for every integer $n$. We denote by $(y'_{n},d''_{n})$
this common value.

Since $M_{l}$ is compact, there are elements $\gamma_{n}$
of $\bar{\Gamma}$ such that ($\gamma_{0}^{n}$ being the linear part of
$\rho(\gamma_{n})$) the $\gamma_{0}^{n}(y_{n},d'_{n})$ converge
to an element $(\bar{y}_{\infty},\bar{d}_{\infty})$ of
$X_{0}$.
Denote by $(\bar{y}',\bar{d}')$ and $(\bar{y}'_{\infty},\bar{d}'_{\infty})$
the images by $\hat{F}$ of $(\bar{y}, \bar{d})$
and $(\bar{y}_{\infty},\bar{d}_{\infty})$.
Then, we have: 

-\qua the $(y'_{n},d''_{n})$ converge to the
element $(\bar{y}',\bar{d}')$ of $X(\rho)$,

-\qua the $\gamma_{n}(y'_{n},d''_{n})$ converge
to the element $(\bar{y}'_{\infty},\bar{d}'_{\infty})$ of $X(\rho)$.

According to Proposition \ref{propac},
the $\gamma_{n}$ are finite in number.
But this is impossible, since the $d'_{n}$ converge
to the projective line $(\bar{y},0)$ and the $\gamma^{n}_{0}d'_{n}$
converge to the projective line $\bar{d}$.
The lemma and the theorem are proven. \fin

\rquea According to Lemma \ref{diff}, we have found a new family of
different differentiable structures on $M_{l}$ for which the
horocyclic foliation remains analytic. The non-triviality of the
moduli of differentiable structures for a given foliation is never an
easy task; this problem has to be compared with the fact that on a
given closed surface there is one and only one differentiable
structure.  Up to our knowledge, the only examples of foliations with
many differentiable structures previously known were the structurally
stable ones, and horocyclic foliations are very far from being
structurally stable!  \erquea

\subsection{Study of the first tautological foliation: the Goldman foliation}
\label{sixdeux}

We prove here the Theorem C.
It is an immediate consequence of Lemmas \ref{neqq}, \ref{empty}
and \ref{whole} below.
Let $\Phi$ be a Goldman foliation on a {\em pure\/} Goldman manifold $M$.
Fix any parametrisation $\Phi^{t}$ of $\Phi$, and any
auxiliary Riemannian metric on $M$.

\begin{defina}
The flow $\Phi^{t}$ is called {\em non-expansive\/} at a point
$x$ of $M$ if, for every $\epsilon>0$, there is
an element $y$ of $M$ and an increasing homeomorphism
$v\co  {\bf R} \rightarrow {\bf R}$ such that:

-\qua $y$ is not on the $\Phi^{t}$--orbit of $x$,

-\qua for any time $t$, the distance between $\Phi^{t}(x)$
and $\Phi^{v(t)}(y)$ is less than $\epsilon$.

The set of points where $\Phi^{t}$ is non-expansive is called the {\em
non-expansiveness locus\/}, and denoted by $\cal N$. Its complement in
$M$ is called the {\em expansiveness locus\/} of $\Phi^{t}$, and
denoted by $\cal E$.  The sets $\cal E$ and $\cal N$ are both
$\Phi^{t}$--invariant.
\end{defina}

\begin{defina}
For any element $(x,d)$ of $X(\rho)$, 
the connected component of $d \cap U(\rho)$ containing $d$
is denoted by $]\alpha(x,d),\beta(x,d)[$. 
\end{defina}

Observe that $\alpha(x,d)$ and $\beta(x,d)$ are
elements of $\Lambda(\rho)$. We choose the notation so
that $\beta(x,d)$ (respectively $\alpha(x,d)$) is the limit when $t$ goes to $+\infty$
(respectively $-\infty$) of $p_{1}(\hat{\Phi}^{t}(x,d))$.

\begin{lema}
\label{neqq}
If $\alpha(x,d) \neq \beta(x,d)$, then the projection $m$
of $(x,d)$ in $M$ belongs to the expansiveness locus.
\end{lema}

\preua
A rigorous and detailed exposition would be long and tedious.
We prefer to indicate the main argument.

Let $K$ be a compact fundamental domain for the action of
$\bar{\Gamma}$ on $X(\rho)$. We consider an ellipse field $\hat{E}$ on
$K$ similar to the ellipse field introduced in the section \ref{foli}:
we fix a euclidean metric on ${\bf R}^{2} \subset {\bf
R}P^{2}_{\ast}$, ie, an ellipse $E_{0}$ and $\hat{E}(x,d)$ is an open
neighborhood of $(x,d)$, on which $p_{2}$ is a homeomorphism with
image the translated of $E_{0}$ centered at $d$.  We extend $\hat{E}$
on the whole $X(\rho)$: for any element $(x,d)$ of $X(\rho)$, the
ellipse $\hat{E}(x,d)$ is $\gamma^{-1}\hat{E}(y,d')$, where $\gamma
(x,d)=(y,d')$ (it can be multidefined for some $(x,d)$, but this has
no incidence for our reasoning).  Define
$E^{t}=p_{2}(\hat{E}(\hat{\Phi}^{t}(x,d)))$: they are ellipses in
${\bf R}^{2}$, centered at $d$. Moreover, they all contain a
subsegment $\sigma_{t}$ centered at $d$ and of length at least
$2\epsilon$ (for the auxiliary euclidean metric on ${\bf R}^{2}$).

For every $t$, let $F_{t}$ be the $p_{2}$--projection of the leaf of
\fr containing $\hat{\Phi}^{t}(x,d)$.  It is a half-plane containing
$E^{t}$, bounded by some line $x(t)$. When $t$ goes towards $+\infty$,
the lines $x(t)$ converge to the line $\beta(x,d)$. Since $x$ belongs
to this limit line, we see that the ellipses $E^{t}$ are more and more
flattened, and converge to the ``degenerated ellipse"
$x^{+}=\beta(x,d)$.

\begin{figure}[htb]
\label{ellipse}
\centerline{\small
%\ShowGrid 
\SetLabels 
\E(0.86*0.9){$x^+$}\\
\E(0.46*0.46){$d$}\\
\endSetLabels 
\AffixLabels{\includegraphics[width=8cm, height=5cm]{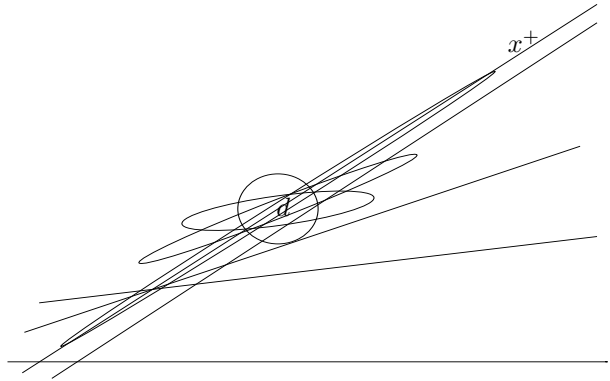}}}
\caption{Ellipses are flattened}
\end{figure}

By the same argument, we see that when $t$ goes to $-\infty$, the
ellipses converge to $x^{-}=\alpha(x,d)$.  When $\alpha(x,d) \neq
\beta(x,d)$, as is assumed in this lemma, the lines $x^{+}$ and
$x^{-}$ are transverse one to the other.  Therefore, the intersection
of all the ellipses $E^{t}$ as $t$ varies over all the real numbers,
the positive and the negative, reduces to $d$.  It follows that
$(x,d)$ belongs to the expansiveness locus: indeed, if the
$\Phi^{t}$--orbit of a point $(x',d')$ remains near the
$\Phi^{t}$--orbit of $(x,d)$, the projection $p_{2}(x',d')=d'$ must
belong to all the ellipses $E^{t}$, and, therefore, $d'$ is equal to
$d$.  \fin

Let $\hat{W}$ be the interior of the set of elements
$(x,d)$ of $X(\rho)$ for which $\alpha(x,d) \neq \beta(x,d)$.
It projects to an open subset $W$ of $M$.
Let $\cal M$ be the complement of $W$
in $M$.

\begin{lema}
\label{empty}
$\hat{W}$ is not empty.
\end{lema}

\preua Let $\gamma_{0}$ be any element of $\bar{\Gamma}$. In ${\bf
R}P^{2}$, it admits $3$ fixed points: $0$, which is of saddle type and
two others which are contained in some projective line $d$ in ${\bf
R}P^{2}$ which is $\gamma_{0}$--invariant. Observe that $d$ meets
$U(\rho)$: if not, $\Lambda(\rho)$ would be contained in $d$, and we
excluded this case while restricting ourselves to pure Goldman
structures.

Let $x$ be any element of $d \cap U(\rho)$.
Then, $(x,d)$ belongs to $X(\rho)$. If $\alpha(x,d) = \beta(x,d)$,
then $d \cap U(\rho)$ is a complete affine line. Therefore,
it must contain at least one fixed point of $\gamma_{0}$, but
this is impossible since the fixed points of $\gamma_{0}$ are outside
$U(\rho)$.

Let $c$ be the subarc of $\Lambda(\rho)$ bounded by 
$\alpha(x,d)$ and $\beta(x,d)$ such that the image of $c$
by the projection of ${\bf R}P^{2} \setminus 0$ along
the leaves of ${\cal G}_{0}$ coincides with the image of
$[\alpha(x,d),\beta(x,d)]$.
The union of
$c$ with $[\alpha(x,d),\beta(x,d)]$ is a Jordan curve
bounding some open subset $V$ of $U(\rho)$ (see figure $3$).
Obviously, the points $(x',d')$ where $x'$ belongs to $V$
and $d'$ is a projective line which doesn't meet $0$ or 
$[\alpha(x,d),\beta(x,d)]$
all satisfy $\alpha(x',d') \neq \beta(x',d')$, and their union is
an open subset of $X(\rho)$. The lemma follows.\fin

\begin{figure}[htb]
\label{piege}\vspace{2mm}
 \centerline{\small
%\ShowGrid 
\SetLabels 
\E(0.1*1){$\Lambda(\rho)$}\\
\E(0.27*0.85){$\alpha(x,d)$}\\
\E(0.31*0.69){$\alpha(x',d')$}\\
\E(0.25*0.35){$\alpha(x'',d'')$}\\
\E(0.8*0.85){$\beta(x,d)$}\\
\E(0.85*0.38){$\beta(x',d')$}\\
\E(0.89*0.57){$\beta(x'',d'')$}\\
\E(0.5*0.85){$d$}\\
\E(0.5*0.6){$d'$}\\
\E(0.41*0.42){$d''$}\\
\E(0.6*0.3){$V$}\\
\E(0.5*0.05){$c$}\\
\endSetLabels 
\AffixLabels{\includegraphics[width=10cm, height=5cm]{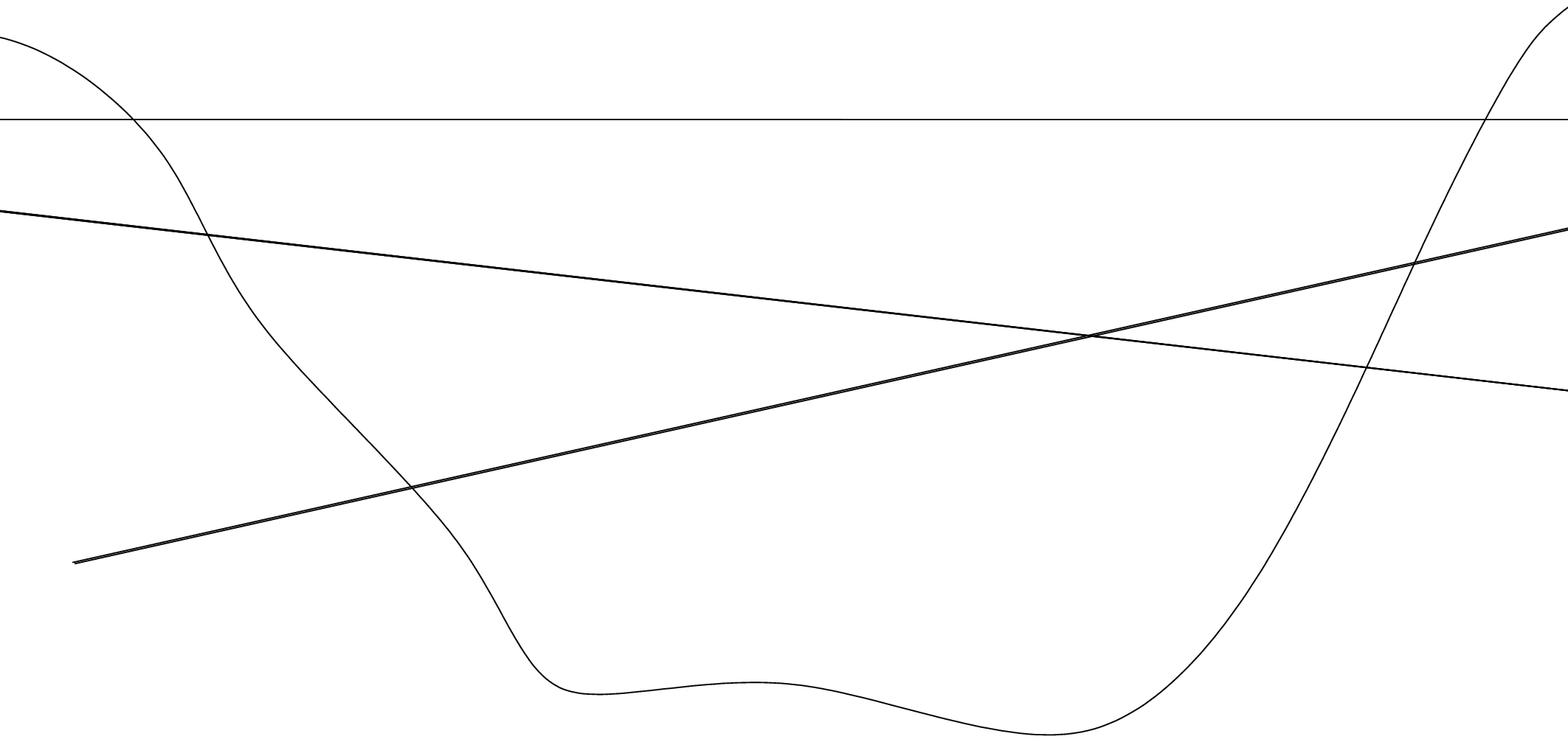}}}
\caption{Exhibiting elements of $W$}
\end{figure}

\begin{lema}
\label{whole}
$\hat{W}$ is not the whole of $X(\rho)$.
\end{lema}

\preua If not, by Lemma \ref{neqq}, the flow $\Phi^{t}$ is expansive.
According to \cite{bru}, it is topologically equivalent to an Anosov
flow. By \cite{ghys}, up to finite coverings, $\Phi^{t}$ is
topologically equivalent to the geodesic flow on the unitary tangent
bundle of a hyperbolic riemmanian surface $S$. There are many ways to
see the impossibility of that. For example: up to finite coverings,
the flow $\hat{\Psi}^{t}$ on $X(\rho)$ must be topologically
equivalent to the geodesic flow of the Poincar\'e disc, and the orbit
space of this geodesic flow is the complement of the diagonal in ${\bf
R}P^{1} \times {\bf R}P^{1}$.  Therefore, the leaf space
$Q_{\hat{\Phi}}$ of $\hat{\Psi}$ is homeomorphic to the annulus, in
particular, it satisfies the Hausdorff separation property.

We observe now that $\Psi^{t}$ is topologically equivalent
to its inverse $\hat{\Psi}^{-t}$ (this property is valid for
any ${\bf R}$--covered Anosov flow whithout cross-section, 
see Theorem C of \cite{bar1}). In particular, for every
perodic orbit, there is another periodic orbit which is freely
homotopic to the initial periodic orbit, but with the inverse
orientation. Let $\gamma_{0}$ be any element of $\bar{\Gamma}$ preserving
an orbit $\theta_{1}$ of $\hat{\Psi}$ (it exists since the geodesic
flow has many periodic orbits). According to the discussion above,
there is another orbit $\theta_{2}$ which is preserved by
$\gamma_{0}$. Let $x_{1}$ and $x_{2}$ be respectively the repelling and
attracting fixed point of $\gamma_{0}$ in ${\bf R}P^{2}$, and let
$d_{0}=(x_{1},x_{2})$ be the projective line containing them.
It is the unique $\gamma_{0}$--invariant projective line which does not contain 
$0$, therefore, it must contain $p_{1}(\theta_{1})$ and $p_{1}(\theta_{2})$.
By $\gamma_{0}$--invariance, it follows that $p_{1}(\theta_{1})$ and $p_{1}(\theta_{2})$
are the two connected components of $d_{0} \setminus \{ x_{1} , x_{2} \}$.
Consider now the intersections of regular leaves of ${\cal G}_{0}$
with $U(\rho)$: we see
them as oriented rays starting from $0$ and reaching $\Lambda(\rho)$.
There are four $\gamma_{0}$--invariant rays, ending at
$x_{1}$ and $x_{2}$. The others rays falls into the two
following exclusive possibilities:

-\qua they meet $p_{1}(\theta_{1})$ or $p_{1}(\theta_{2})$,

-\qua they meet $\Lambda(\rho)$ before meeting $p_{1}(\theta_{1})$ 
or $p_{1}(\theta_{2})$.

It is easy to see that the union of the rays of the first category
contains the union of two triangles with vertices $0$, $x_{1}$ and
$x_{2}$. In other words, there is an affine half-plane $T$ contained
in $U(\rho)$ bounded by $d_{0}$ and another projective line containing
$0$ and $x_{i}$ (where $i=1$ or $2$).  Let $d$ be any projective line
containing $x_{i}$ and intersecting $T$: the affine line $d \setminus
\{ x_{i} \}$ is the projection by $p_{1}$ of some leaf of
$\hat{\Psi}$. This line $d$ can be arbitrarly close to $d_{0}$, and
thus, arbitrarly close to $p_{1}(\theta_{1})$ and $p_{1}(\theta_{2})$.
This is in contradiction with the Hausdorff separation property in
$Q_{\hat{\Psi}}$.  \fin

\rquea
The proof of \ref{whole} we propose here uses very deep results.
We don't know a more elementary one.
\erquea

\section{Conclusion}
Pure Goldman foliations are good examples of analytic foliations,
satisfying strong properties, but which remain quite mysterious.  Many
questions remain open.  It would be worthwhile to know a little more
about them.

{\bf Question 1}\qua Does a Goldman foliation admit periodic
orbits?  Observe that such a periodic orbit is necessarily
hyperbolic. Observe also that if a pure Goldman foliation has no
periodic orbits and every element of the affine holonomy group is of
determinant $1$, it would be a non-minimal analytical volume
preserving foliation without periodic orbit on a $3$--manifold.  Such
foliations are not so easy to construct, the first known example being
the Kuperberg foliation \cite{kuper}.

{\bf Question 2}\qua We proved that a pure Goldman foliation is not
minimal by exhibiting a non-trivial closed invariant subset $\cal M$.
Is $\cal M$ itself minimal? How can we describe the dynamic of the
Goldman foliation on $\cal M$?

{\bf Question 3}\qua Being conjugate to horocyclic foliations, the
first tautological foliations of pure Goldman flag manifolds are
uniquely ergodic: there is a unique invariant measure.  When is this
measure absolutely continuous with respect to the Lebesgue measure?

{\bf Question 4}\qua What are the ergodic properties of
Goldman foliations?  

{\bf Question 5}\qua We can suspect, from the expansiveness of a
pure Goldman flow outside $\cal M$, that its measure entropy is
positive. Is this true? This question is related to the question $1$,
since, according to a theorem by A Katok, the entropy of a regular
flow on a closed $3$--manifold without periodic orbit is zero
\cite{katok}.  The positivity of the entropy would follow if we could
show that the Lyapounov exponents are not all zero almost
everywhere. The paper \cite{guivarch} of Y Guivarc'h establishes some
results in this direction.  Unfortunately, they apply to groups of
projective transformations which do not preserve any projective
subspaces, which is certainly not the case for the groups we have
considered here.

{\bf Question 6}\qua We know that the Jordan curve $\Lambda(\rho)$
is not Lipschitz. But we can wonder what is its regularity. Is it
H\"{o}lder?  Is it rectifiable?

{\bf Question 7}\qua In Theorem B, can we withdraw the assumption
forcing the image of the developing image to be contained in
$X_{\infty}$?  In other words, is it true that any flag structure on a
Seifert manifold, for which the holonomy morphism is hyperbolic, is a
finite covering of $M_{H}$? The answer is expected to be yes.

{\bf Question 8}\qua We only considered deformations of holonomy
groups inside $Af_{0}$. What happens for general deformations inside
the whole $SL(3,{\bf R})$?  Do they still act freely and properly
discontinuously on some open subset of $X$ with compact quotient?

\end{document}

%% file: gtoutput.tex
\def\ifplaintex{\expandafter\ifx\csname documentclass\endcsname\relax}

\ifplaintex 
\else
\fi

\expandafter\ifx\csname beginpicture\endcsname\relax
\expandafter\ifx\csname documentclass\endcsname\relax
\input pictex \else
\input prepictex \input pictex \input postpictex \fi\fi

\def\gt{{\mathsurround=0pt\it $\cal G\mskip-2mu$eometry \&\ 
$\cal T\!\!$opology}}        %  journal title in recommended style

\def\gtp{{\mathsurround=0pt\it $\cal G\mskip-2mu$eometry \&\ 
$\cal T\!\!$opology $\cal P\!$ublications}}  % GT publications

\def\volumenumber#1{\def\thevolumenumber{#1}}
\def\papernumber#1{\def\thepapernumber{#1}}
\def\volumeyear#1{\def\thevolumeyear{#1}}

\def\pagenumbers#1#2{\def\startpage{#1}\def\finishpage{#2}}
\def\published#1{\def\publishdate{#1}}
\def\proposed#1{\def\theproposer{#1}}
\def\seconded#1{\def\theseconders{#1}}
\def\received#1{\def\receiveddate{#1}}
\def\revised#1{\def\reviseddate{#1}}
\def\accepted#1{\def\accepteddate{#1}}

\def\asciiaddress#1{\def\theasciiaddress{#1}}

\long\def\asciiabstract#1{\long\def\theasciiabstract{#1}}

\let\\\par
\let\thevolumenumber\relax\let\thepapernumber\relax
\let\thevolumeyear\relax\let\thesamplenumber\relax\let\startpage\relax
\let\finishpage\relax\let\publishdate\relax\let\receiveddate\relax
\let\reviseddate\relax\let\accepteddate\relax\let\theasciititle\relax
\let\theasciiauthors\relax\let\theasciiaddress\relax
\let\theasciiabstract\relax
\let\theasciiemail\relax\let\theshortauthors\relax\let\theshorttitle\relax

\long\def\maketitlep{   % start of definition of \maketitlep

\count0=\startpage

\gt\hfill      %   Journal title (top left) 
\beginpicture
\setcoordinatesystem units <0.33truein, 0.33truein> point at 2.2 0.9
\setplotsymbol ({$\cal G$})
\plotsymbolspacing=9truept
\circulararc 315 degrees from 0 1 center at 0 0
\setplotsymbol ({$\cal T$})
\circulararc 315 degrees from 1 -1 center at 1 0
\endpicture
\break
{\small\ifx\thesamplenumber\relax % sample?  
Volume \else Sample
\fi\thevolumenumber\ (\thevolumeyear)
\startpage--\finishpage\nl
Published: \publishdate}
\vglue 0.5truein plus 0.4fil minus 0.1truein

{\parskip=0pt\leftskip 0pt plus 1fil\def\\{\par\smallskip}{\ifplaintex\large
\else\Large\fi\bf\thetitle}\par\medskip}   

\vglue 0pt plus 0.1fil 

{\parskip=0pt\leftskip 0pt plus 1fil\def\\{\par}{\sc\theauthors}
\par\medskip}

\vglue 0pt plus 0.1fil 

{\small\parskip=0pt\let\newline\\
{\leftskip 0pt plus 1fil\def\\{\par}{\sl\theaddress}\par}
\expandafter\ifx\theemail\relax    % email address?
\relax\else\vglue 5pt plus 0.02fil minus 2pt\def\\{\stdspace{\rm 
and}\stdspace} 
\cl{Email:\stdspace\tt\theemail}\fi
\ifx\theurl\relax                  % URL given?
\relax\else\vglue 5pt plus 0.02fil minus 2pt\def\\{\stdspace{\rm 
and}\stdspace}
\cl{URL:\stdspace\tt\theurl}\fi\par}

\vglue 7pt plus 0.3fil minus 3pt

{\bf Abstract}
\vglue 5pt plus 0.1fil minus 2pt

\theabstract

\vglue 7pt plus 0.3fil minus 3pt

{\bf AMS Classification numbers}\quad Primary:\quad \theprimaryclass

Secondary:\quad \thesecondaryclass

\vglue 5pt plus 0.3fil minus 2pt

{\bf Keywords}\quad \thekeywords

\vglue 10pt plus 0.5fil minus 5pt

{\small  Proposed: \theproposer\hfill Received: \receiveddate\nl
Seconded: \theseconders\hfill 
\ifx\reviseddate\relax                         % paper revised?
Accepted: \accepteddate                        % no
\else
Revised: \reviseddate                          % yes
\fi}
\eject
}       %  end of definition of \maketitlep

\let\maketitlepage\maketitlep
\let\maketitle\maketitlepage

\font\phead=cmsl9 scaled 950
\font=cmsl9 scaled 1050
\font\pnum=cmbx10 scaled 913
\font\lnum=cmbx10 
\font\pfoot=cmsl9 scaled 950
\font=cmsl9 scaled 1050
\ifplaintex
\headline{\vbox to 0pt{\vskip -4.5mm\line{\small\phead\ifnum
\count0=\startpage ISSN 1364-0380 (on line)
1465-3060 (printed) \hfill {\pnum\folio}\else\ifodd\count0\def\\{ }% 
\ifx\theshorttitle\relax\thetitle\else\theshorttitle\fi\hfill{\pnum\folio}
\else\def\\{ and }{\pnum\folio}\hfill\ifx\theshortauthors\relax\theauthors
\else\theshortauthors\fi\fi\fi}\vss}}
\footline{\vbox to 0pt{\vglue 0mm\line{\small\pfoot\ifnum\count0=\startpage
\copyright\ \gtp\hfill\else
\gt, Volume \thevolumenumber\ (\thevolumeyear)\hfill\fi}\vss
}}
\else
\makeatletter
\def\@oddhead{{\small\lhead\ifnum\count0=\startpage ISSN 1364-0380 (on line)
1465-3060 (printed) \hfill {\lnum\number\count0}\else\ifodd\count0
\def\\{ }\ifx\theshorttitle\relax \thetitle \else\theshorttitle\fi\hfill
{\lnum\number\count0}\else\def\\{ and }{\lnum\number\count0}
\hfill\ifx\theshortauthors\relax 
\theauthors\else\theshortauthors\fi\fi\fi}}\def\@evenhead{\@oddhead}
\def\@oddfoot{\small\lfoot\ifnum\count0=\startpage\copyright\ \gtp\hfill\else
\gt, Volume \thevolumenumber\ (\thevolumeyear)\hfill\fi}
\def\@evenfoot{\@oddfoot}
\makeatother
\fi

\newwrite\gtoutfile
\long\gdef\makeheadfile{  %%% start of definition of \makeheadfile
{\def\\{, }\def\s{ }
\immediate\openout\gtoutfile head.xxx
\immediate\write\gtoutfile{To: math@arxiv.org}
\immediate\write\gtoutfile{Subject: put or rep NNNNN:pppp}
\immediate\write\gtoutfile{--text follows this line--}
\immediate\write\gtoutfile{Proxy-for: \ifx\theasciiauthors\relax
\theauthors\else\theasciiauthors\fi\s<\ifx\theasciiemail\relax\theemail\else\theasciiemail\fi>}
\immediate\write\gtoutfile{\noexpand\\}
\immediate\write\gtoutfile{Authors: \ifx\theasciiauthors\relax
\theauthors\else\theasciiauthors\fi}
\immediate\write\gtoutfile{Title: \ifx\theasciititle\relax
\thetitle\else\theasciititle\fi}
\immediate\write\gtoutfile{Subj-class: GT or SG or MG etc}
\immediate\write\gtoutfile{MSC-class: \theprimaryclass\ifx\thesecondaryclass\relax\else, \thesecondaryclass\fi}
\immediate\write\gtoutfile{Journal-ref: Geom. Topol. \thevolumenumber
(\thevolumeyear) \startpage-\finishpage}
\immediate\write\gtoutfile{Comments: Published by Geometry and Topology at}
\immediate\write\gtoutfile{\s\s http://www.maths.warwick.ac.uk/gt/GTVol\thevolumenumber/paper\thepapernumber.abs.html}
\immediate\write\gtoutfile{\noexpand\\}
\immediate\write\gtoutfile{}
\ifx\theasciiabstract\relax
\immediate\write\gtoutfile{\theabstract}\else
\immediate\write\gtoutfile{\theasciiabstract}\fi
\immediate\write\gtoutfile{}
\immediate\write\gtoutfile{\noexpand\\}
\immediate\write\gtoutfile{}
\immediate\closeout\gtoutfile}}  %%% end of definition of \makeheadfile

\def\maketitlepage{\maketitlep\makeheadfile}
\let\maketitle\maketitlepage

\def\ifplaintex{\expandafter\ifx\csname documentclass\endcsname\relax}

\ifplaintex 
\else
\fi

\expandafter\ifx\csname beginpicture\endcsname\relax
\expandafter\ifx\csname documentclass\endcsname\relax
\input pictex \else
\input prepictex \input pictex \input postpictex \fi\fi

\def\gt{{\mathsurround=0pt\it $\cal G\mskip-2mu$eometry \&\ 
$\cal T\!\!$opology}}        %  journal title in recommended style

\def\gtp{{\mathsurround=0pt\it $\cal G\mskip-2mu$eometry \&\ 
$\cal T\!\!$opology $\cal P\!$ublications}}  % GT publications

\def\volumenumber#1{\def\thevolumenumber{#1}}
\def\papernumber#1{\def\thepapernumber{#1}}
\def\volumeyear#1{\def\thevolumeyear{#1}}

\def\pagenumbers#1#2{\def\startpage{#1}\def\finishpage{#2}}
\def\published#1{\def\publishdate{#1}}
\def\proposed#1{\def\theproposer{#1}}
\def\seconded#1{\def\theseconders{#1}}
\def\received#1{\def\receiveddate{#1}}
\def\revised#1{\def\reviseddate{#1}}
\def\accepted#1{\def\accepteddate{#1}}

\def\asciiaddress#1{\def\theasciiaddress{#1}}

\long\def\asciiabstract#1{\long\def\theasciiabstract{#1}}

\let\\\par
\let\thevolumenumber\relax\let\thepapernumber\relax
\let\thevolumeyear\relax\let\thesamplenumber\relax\let\startpage\relax
\let\finishpage\relax\let\publishdate\relax\let\receiveddate\relax
\let\reviseddate\relax\let\accepteddate\relax\let\theasciititle\relax
\let\theasciiauthors\relax\let\theasciiaddress\relax
\let\theasciiabstract\relax
\let\theasciiemail\relax\let\theshortauthors\relax\let\theshorttitle\relax

\long\def\maketitlep{   % start of definition of \maketitlep

\count0=\startpage

\gt\hfill      %   Journal title (top left) 
\beginpicture
\setcoordinatesystem units <0.33truein, 0.33truein> point at 2.2 0.9
\setplotsymbol ({$\cal G$})
\plotsymbolspacing=9truept
\circulararc 315 degrees from 0 1 center at 0 0
\setplotsymbol ({$\cal T$})
\circulararc 315 degrees from 1 -1 center at 1 0
\endpicture
\break
{\small\ifx\thesamplenumber\relax % sample?  
Volume \else Sample
\fi\thevolumenumber\ (\thevolumeyear)
\startpage--\finishpage\nl
Published: \publishdate}
\vglue 0.5truein plus 0.4fil minus 0.1truein

{\parskip=0pt\leftskip 0pt plus 1fil\def\\{\par\smallskip}{\ifplaintex\large
\else\Large\fi\bf\thetitle}\par\medskip}   

\vglue 0pt plus 0.1fil 

{\parskip=0pt\leftskip 0pt plus 1fil\def\\{\par}{\sc\theauthors}
\par\medskip}

\vglue 0pt plus 0.1fil 

{\small\parskip=0pt\let\newline\\
{\leftskip 0pt plus 1fil\def\\{\par}{\sl\theaddress}\par}
\expandafter\ifx\theemail\relax    % email address?
\relax\else\vglue 5pt plus 0.02fil minus 2pt\def\\{\stdspace{\rm 
and}\stdspace} 
\cl{Email:\stdspace\tt\theemail}\fi
\ifx\theurl\relax                  % URL given?
\relax\else\vglue 5pt plus 0.02fil minus 2pt\def\\{\stdspace{\rm 
and}\stdspace}
\cl{URL:\stdspace\tt\theurl}\fi\par}

\vglue 7pt plus 0.3fil minus 3pt

{\bf Abstract}
\vglue 5pt plus 0.1fil minus 2pt

\theabstract

\vglue 7pt plus 0.3fil minus 3pt

{\bf AMS Classification numbers}\quad Primary:\quad \theprimaryclass

Secondary:\quad \thesecondaryclass

\vglue 5pt plus 0.3fil minus 2pt

{\bf Keywords}\quad \thekeywords

\vglue 10pt plus 0.5fil minus 5pt

{\small  Proposed: \theproposer\hfill Received: \receiveddate\nl
Seconded: \theseconders\hfill 
\ifx\reviseddate\relax                         % paper revised?
Accepted: \accepteddate                        % no
\else
Revised: \reviseddate                          % yes
\fi}
\eject
}       %  end of definition of \maketitlep

\let\maketitlepage\maketitlep
\let\maketitle\maketitlepage

\font\phead=cmsl9 scaled 950
\font=cmsl9 scaled 1050
\font\pnum=cmbx10 scaled 913
\font\lnum=cmbx10 
\font\pfoot=cmsl9 scaled 950
\font=cmsl9 scaled 1050
\ifplaintex
\headline{\vbox to 0pt{\vskip -4.5mm\line{\small\phead\ifnum
\count0=\startpage ISSN 1364-0380 (on line)
1465-3060 (printed) \hfill {\pnum\folio}\else\ifodd\count0\def\\{ }% 
\ifx\theshorttitle\relax\thetitle\else\theshorttitle\fi\hfill{\pnum\folio}
\else\def\\{ and }{\pnum\folio}\hfill\ifx\theshortauthors\relax\theauthors
\else\theshortauthors\fi\fi\fi}\vss}}
\footline{\vbox to 0pt{\vglue 0mm\line{\small\pfoot\ifnum\count0=\startpage
\copyright\ \gtp\hfill\else
\gt, Volume \thevolumenumber\ (\thevolumeyear)\hfill\fi}\vss
}}
\else
\makeatletter
\def\@oddhead{{\small\lhead\ifnum\count0=\startpage ISSN 1364-0380 (on line)
1465-3060 (printed) \hfill {\lnum\number\count0}\else\ifodd\count0
\def\\{ }\ifx\theshorttitle\relax \thetitle \else\theshorttitle\fi\hfill
{\lnum\number\count0}\else\def\\{ and }{\lnum\number\count0}
\hfill\ifx\theshortauthors\relax 
\theauthors\else\theshortauthors\fi\fi\fi}}\def\@evenhead{\@oddhead}
\def\@oddfoot{\small\lfoot\ifnum\count0=\startpage\copyright\ \gtp\hfill\else
\gt, Volume \thevolumenumber\ (\thevolumeyear)\hfill\fi}
\def\@evenfoot{\@oddfoot}
\makeatother
\fi

\newwrite\gtoutfile
\long\gdef\makeheadfile{  %%% start of definition of \makeheadfile
{\def\\{, }\def\s{ }
\immediate\openout\gtoutfile head.xxx
\immediate\write\gtoutfile{To: math@arxiv.org}
\immediate\write\gtoutfile{Subject: put or rep NNNNN:pppp}
\immediate\write\gtoutfile{--text follows this line--}
\immediate\write\gtoutfile{Proxy-for: \ifx\theasciiauthors\relax
\theauthors\else\theasciiauthors\fi\s<\ifx\theasciiemail\relax\theemail\else\theasciiemail\fi>}
\immediate\write\gtoutfile{\noexpand\\}
\immediate\write\gtoutfile{Authors: \ifx\theasciiauthors\relax
\theauthors\else\theasciiauthors\fi}
\immediate\write\gtoutfile{Title: \ifx\theasciititle\relax
\thetitle\else\theasciititle\fi}
\immediate\write\gtoutfile{Subj-class: GT or SG or MG etc}
\immediate\write\gtoutfile{MSC-class: \theprimaryclass\ifx\thesecondaryclass\relax\else, \thesecondaryclass\fi}
\immediate\write\gtoutfile{Journal-ref: Geom. Topol. \thevolumenumber
(\thevolumeyear) \startpage-\finishpage}
\immediate\write\gtoutfile{Comments: Published by Geometry and Topology at}
\immediate\write\gtoutfile{\s\s http://www.maths.warwick.ac.uk/gt/GTVol\thevolumenumber/paper\thepapernumber.abs.html}
\immediate\write\gtoutfile{\noexpand\\}
\immediate\write\gtoutfile{}
\ifx\theasciiabstract\relax
\immediate\write\gtoutfile{\theabstract}\else
\immediate\write\gtoutfile{\theasciiabstract}\fi
\immediate\write\gtoutfile{}
\immediate\write\gtoutfile{\noexpand\\}
\immediate\write\gtoutfile{}
\immediate\closeout\gtoutfile}}  %%% end of definition of \makeheadfile

\def\maketitlepage{\maketitlep\makeheadfile}
\let\maketitle\maketitlepage

%% file: 2001-7.bbl
\newpage

\begin{thebibliography}

\bibitem{abe} {\bf Ryuji Abe}, {\em Geometric approach to rigidity of
horocyclics\/}, Tokyo J. Math. {18} (1995) 271--283

\bibitem{norm} {\bf V Bangert}, {\em Minimal geodesics\/}, Ergod. Th. {\&}
Dynam. Sys. {10} (1989) 263--286

\bibitem{barbot} {\bf T Barbot}, {\em Actions de groupes sur les
1--vari\'et\'es non s\'epar\'ees et feuilletages de codimension un\/},
Ann. Fac. Sci. Toulouse, {7} (1998) 559--597

\bibitem{bar1} {\bf T. Barbot}, {\em Caract\'erisation des flots
d'Anosov en dimension 3 par leurs feuilletages faibles\/},
Ergod. Th. {\&} Dynam. Sys. {15} (1995) 247--270 

\bibitem{benzecri} {\bf J\,P Benz\'ecri}, {\em Vari\'et\'es localement
affines et projectives\/}, Bull. Soc. Math. France {88} (1960)
229--332

\bibitem{bers} {\bf L Bers}, {\em Spaces of Kleinian groups\/},
Lectures Notes in Math. {155} (1970) 9--34

\bibitem{bowen} {\bf R Bowen}, {\bf B Marcus}, {\em Unique ergodicity
for horocycle foliations\/}, Israel J. Math. {26} (1977)
43--67

\bibitem{bru}{\bf M Brunella}, {\em Expansive foliations on Seifert
manifolds and on torus bundles\/}, Bol. Soc. Brasil. Math. {24} (1993)
89--104

\bibitem{canary} {\bf R\,D Canary}, {\bf D\,B\,A Epstein}, {\bf P
Green}, {\em Notes on notes of Thurston\/}, from: ``Analitical and
Geometric Aspects of Hyperbolic Space'', (D\,B\,A Epstein, editor)
London Math. Soc. Lect. Notes Series {111} (1986) 3--92

\bibitem{carriere} {\bf Y Carri\`ere}, {\em Autour de la conjecture de
L Markus sur les vari\'et\'es affines\/}, Invent. Math. {95} (1989)
615--628

\bibitem{cdcr1} {\bf S Choi}, {\em Convex decompositions of real
projective surfaces}. I: {\em $\pi$--annuli and convexity},
J. Differential Geom. {40} (1994) 165--208

\bibitem{cdcr2} {\bf S Choi}, {\em Convex decompositions of real projective
surfaces}. II: {\em Admissible decompositions},
J. Differential Geom. {40} (1994) 239--283 

\bibitem{cdcr3} {\bf S Choi}, {\em Convex decompositions of real
projective surfaces}. III: {\em For closed and nonorientable
surfaces}, J. Korean Math. Soc. {33} (1996) 1138--1171

\bibitem{CG} {\bf S Choi}, {\bf W\,M Goldman}, {\em The classification
of real projective structures on compact surfaces},
Bull. Amer. Math. Soc.  {34} (1997) 161--171

\bibitem{choigold} {\bf S Choi}, {\bf W\,M Goldman}, {\em Convex real
structures on closed surfaces are closed\/},
Proc. Amer. Math. Soc. {118} (1993) 657--661

\bibitem{ghys} {\bf E Ghys}, {\em Flots d'Anosov sur les
3--vari\'et\'es fibr\'ees en cercles\/}, Ergod. Th. {\&}
Dynam. Sys. {4} (1984) 67--80

\bibitem{ghyqf2} {\bf E Ghys}, {\em Flots d'Anosov dont les
feuilletages stables sont différentiables\/}, Ann. Sci. \'Ecole
Norm. Sup. {20} (1987) 251--270

\bibitem{golbou} {\bf W Goldman}, {\em Geometric structures on
manifolds and varieties of representations\/}, Contemp. Math. {74}
(1988) Amer. Math. Soc. Providence, RI

\bibitem{guivarch} {\bf Y Guivarc'h}, {\em Produits de matrices
al\'eatoires et applications aux propri\'et\'es g\'eom\'etriques des
sous-groupes du groupe lin\'eaire\/}, Ergod. Th. {\&} Dynam. Sys. {10}
(1990) 483--512

\bibitem{haefliger} {\bf A Haefliger}, {\em Groupo\"{\i}des
d'holonomie et classifiants\/}, Ast\'erisque {116} (1984) 70--97

\bibitem{katok} {\bf A Katok}, {\em Lyapounov exponents, entropy and
periodic orbits for diffeomorphisms\/}, Inst. Hautes \'Etudes
Sci. Publ. Math. {51} (1980) 137--173

\bibitem{kosz} {\bf J\,L Koszul}, {\em Déformations de connexions
localement plates\/}, Ann. Inst. Fourier (Grenoble) {18} (1968) fasc. 1,
103--114

\bibitem{kuper} {\bf G Kuperberg}, {\bf K Kuperberg}, {\em Generalized
counterexamples to the Seifert conjecture\/}, Ann. of Math. {144}
(1996) 239--268

\bibitem{marden} {\bf A Marden}, {\em The geometry of finitely
generated kleinian groups\/}, Ann. of Math. {99} (1974) 383--462

\bibitem{matsu} {\bf S Matsumoto}, {\em Affine foliations on
3--manifolds\/}, preprint, Nihon University

\bibitem{naganoyagi} {\bf T Nagano}, {\bf K Yagi}, {\em The affine
structures on the real two-torus\/}, Osaka J. Math. {11} (1974)
181--210

\bibitem{plante3} {\bf J\,F Plante}, {\em Anosov flows\/}, 
Amer. J. Math. {94} (1972) 729--754


\bibitem{ratcliffe} {\bf J\,G Ratcliffe}, {\em Foundations of
hyperbolic manifolds\/}, Grad. Texts in Math. {149}, Springer--Verlag


\bibitem{salein} {\bf F Salein}, {\em Vari\'et\'es anti-de Sitter de
dimension 3 poss\'edant un champ de Killing non trivial\/},
C. R. Acad. Sci. Paris S\`er. I Math. {324} (1997) 525--530

\bibitem{sulli} {\bf D Sullivan}, cours I.H.E.S.

\bibitem{thuthese} {\bf W Thurston}, {\em Foliations on 3--manifolds
which are circle bundles\/}, Thesis, Berkeley (1972).


\end{thebibliography}
